\documentclass[10pt]{amsart}
\usepackage{amscd,amssymb,hyperref}
\usepackage[PostScript=dvips]{diagrams}

\usepackage{graphics}
\usepackage{epsfig}
\usepackage{picinpar}
\usepackage{wrapfig}
\usepackage{amsmath}
\usepackage{color}
\usepackage{pstricks,pst-node,pst-text,pst-3d}
\usepackage{pst-plot}
\usepackage{verbatim}

\newtheorem{thm}{Theorem}[section]
\newtheorem{lem}[thm]{Lemma}
\newtheorem{cor}[thm]{Corollary}
\newtheorem{prop}[thm]{Proposition}

\newtheorem{rem}{Remark}[section]
\newtheorem{example}{Example}[section]

\def\square{\vbox{
      \hrule height 0.4pt
      \hbox{\vrule width 0.4pt height 5.5pt \kern 5.5pt \vrule width 0.4pt}
      \hrule height 0.4pt}}

\def\colim{\operatorname{co l i m}}
\def\id{\mathrm{id}}

\def\Ker{\mathrm{K er}}
\def\ch\mathrm{c h}

\def\hocolim{\operatorname{ho co l i m}}

\newcommand{\Z}{\mathbb{Z}}

\newcommand{\R}{\ensuremath{\mathbb{R}}}

\newcommand{\calC}{\ensuremath{\mathcal{C}}}

\newcommand{\calB}{\ensuremath{\mathcal{B}}}

\newcommand{\calZ}{\ensuremath{\mathcal{Z}}}

\newcommand{\calA}{\ensuremath{\mathcal{A}}}

\newcommand{\calG}{\ensuremath{\mathcal{G}}}

\let\la=\langle
\let\ra=\rangle

\numberwithin{equation}{section}

\begin{document}

\newcommand{\auths}[1]{\textrm{#1},}
\newcommand{\artTitle}[1]{\textsl{#1},}
\newcommand{\jTitle}[1]{\textrm{#1}}
\newcommand{\Vol}[1]{\textbf{#1}}
\newcommand{\Year}[1]{\textrm{(#1)}}
\newcommand{\Pages}[1]{\textrm{#1}}

\author{J. Wu$^{\dag}$}
\address{Department of Mathematics\\
National University of Singapore\\
Singapore 119260\\
Republic of Singapore}\email{matwujie@math.nus.edu.sg}
\urladdr{http://www.math.nus.edu.sg/\~{}matwujie}
\thanks{$^{\dag}$ Research is supported in part by the Academic Research Fund of the
National University of Singapore R-146-000-101-112.}


\begin{abstract}
We introduce the (general) homotopy groups of spheres as link invariants for Brunnian-type links through the investigations on the intersection subgroup of the normal closures of the meridians of strongly nonsplittable links. The homotopy groups measure the difference between the intersection subgroup and symmetric commutator subgroup of the normal closures of the meridians and give the invariants of the links obtained in this way. Moreover the higher homotopy-group invariants can produce some links that could not be detected by the Milnor invariants. Furthermore all homotopy groups of spheres can be obtained from the geometric Massey products on links.
\end{abstract}

\title[Brunnian-type links and homotopy groups of spheres]{On Brunnian-type links and the link invariants given by homotopy groups of spheres}

\maketitle

\section{Introduction}
Let $L=\{l_1,l_2,\ldots,l_n\}$ be an $n$-link in $S^3$, where $l_i$ is the $i\,$th component of $L$. The link group $G(L)$ is defined to be the fundamental group of the link complement $S^3\smallsetminus L$. Let
$$
\{l_1,l_2,\ldots,l_{i-1},l_{i+1},\ldots,l_n\}
$$
be the $(n-1)$-link obtained by removing the $i\,$th link component of $L$, denoted by $d_iL$ or $d_{l_i}L$. The inclusion of link complement
$
S^3\smallsetminus L\hookrightarrow S^3 \smallsetminus d_iL
$
induces a group homomorphism
$$
d_i\colon G(L)\longrightarrow G(d_iL).
$$
Let $A(L,l_i)$ be the kernel of the homomorphism $d_i$. Note that $A(L,l_i)$ is the normal closure of the meridian of $l_i$ in $G(L)$. Let $L'=\{l_{i_1},l_{i_2},\ldots,l_{i_t}\}$ be a sublink of $L$. Consider the intersection subgroup
\begin{equation}\label{equation1.1}
A(L,L')=\bigcap_{j=1}^tA(L,l_{i_j}).
\end{equation}
Given any element $\alpha \in A(L,L')$, one can choose a knot $K$ in the link component $S^3\smallsetminus L$ as a representative for the homotopy class $\alpha$. The $(n+1)$-link $\tilde L=L\cup K$ admits the \textit{Brunnian-type property} that the knot $K$ becomes a trivial knot up to pointed homotopy in the link complement $S^3-d_{i_j}L$ for $1\leq j\leq t$. Namely, by removing any $i_j\,$th component of $L$, the knot $K$ separates away from $d_iL$ and is homotopic to the trivial knot up to pointed homotopy. Thus the intersection subgroup~(\ref{equation1.1}) helps to construct Brunnian-type links.

There is a canonical subgroup of the intersection subgroup~(\ref{equation1.1}) given by iterated commutators. Let
$$
A_S[L,L']=\prod_{\sigma\in \Sigma_t}[[A(L,l_{i_{\sigma(1)}}),A(L,l_{i_{\sigma(2)}})],\ldots,A(L,l_{i_{\sigma(t)}})]
$$
be the symmetric (iterated) commutator subgroup of the normal subgroups of $A(L,l_{i_j})$ for $1\leq j\leq t$. Note that the group $A_S[L,L']$ is generated by the $t$-fold iterated commutators
$$
[[g_1,g_2],\ldots,g_t]
$$
with $g_j\in A(L,l_{i_{\sigma(j)}})$, $1\leq j\leq t$, for some permutation $\sigma\in\Sigma_t$. Clearly the symmetric commutator subgroup
$A_S[L,L']$ is a (normal) subgroup of $\bigcap_{j=1}^tA(L,l_{i_j})$
because
$$
d_{i_j}([[g_1,g_2],\ldots,g_t])=[[d_{i_j}(g_1),d_{i_j}(g_2)],\ldots,d_{i_j}(g_t)]=1
$$
for each $1\leq j\leq t$.

The purpose of this article is to investigate the quotient group
\begin{equation}\label{equation1.2}
\calA(L,L')=\frac{A(L,L')}{A_S[L,L']}.
\end{equation}
Observe that the intersection subgroup ~(\ref{equation1.1}) is given by the short exact sequence
$$
A_S[L,L']\rightarrowtail A(L,L')\twoheadrightarrow\calA(L,L').
$$
with a set of generators for $A_S[L,L']$ being understood in some sense. Our determination of the group $\calA(L,L')$ will be given in terms of the homotopy groups of spheres for some pairs of links $(L,L')$.

Recall that a link $L$ is \textit{splittable} if there is an embedding of $S^2$ into $S^3$ such that each side of $S^3\smallsetminus S^2$ has nontrivial intersection of $L$. We call a link $L$ \textit{strongly nonsplittable} if any nonempty sublink of $L$ is not splittable. Observe that any link $L$ has a decomposition as a disjoint union of nonsplittable sublinks. Let $\nu+1$ be the number of nonsplittable components of $L$. We call the number $\nu$ the \textit{splitting genus} of $L$ denoted by $\nu(L)$ with $\nu(L)\geq 0$. A pair of links $(L,L_0)$ in $S^3$ means a link $L$ with a sublink $L_0$. We call $(L,L_0)$ \textit{strongly nonsplittable} if any sublink $L'$ of $L$ with $L_0\subsetneq L'$ is not splittable. Intuitively the link $L$ is obtained from $L_0$ by adding link components with each of them going around all nonsplittable components of $L_0$. For instance let $L_0$ be a trivial 2-link. We may regard $L_0$ as the two separate rings holding on a tree. Let $L$ be any $(n+2)$-link by adding $n$ link components with each of them going around the both rings on the tree. Then $(L,L_0)$ is strongly nonsplittable. The nontrivial elements in the intersection subgroup $A(L,L\smallsetminus L_0)$ can be intuitively described as an extra loop going around the link $L$ which would fall down to the ground by removing any one of the link components in $L\smallsetminus L_0$. Our main result is as follows:

\begin{thm}\label{theorem1.1}
Let $(L,L_0)$ be a pair of links in $S^3$ such that $L\smallsetminus L_0$ is an $n$-link with $n\geq 2$ and $(L,L_0)$ is strongly nonsplittable. Let $L'$ be any sub $t$-link of $L\smallsetminus L_0$ with $2\leq t\leq n$.
\begin{enumerate}
\item If $L'\subsetneq L\smallsetminus L_0$ or $L_0$ is nonempty and nonsplittable, then $\calA(L,L')=0$.
\item If $L_0=\emptyset$ and $L'=L$, then there is an isomorphism of groups
$$
\calA(L,L)\cong \pi_n(S^3).
 $$
 Moreover any element $\alpha \in \calA(L,L)\cong \pi_n(S^3)$ gives a link invariant for the $(n+1)$-links $L\cup l$ with $l$ a knot in $S^3\smallsetminus L$ representing $\alpha$.
\item If $L'=L\smallsetminus L_0$ and $L_0$ is nonempty and splittable with splitting genus $\nu\geq1$, then there is an isomorphism of $\Z(G(L))$-modules
$$
\calA(L,L\smallsetminus L_0)\cong \pi_n\left( \bigvee_{j=1}^\nu G(L_0)\ltimes S^2\right)\cong \pi_n\left(\bigvee_{j=1}^\nu S^2\right)\oplus\pi_n\left( \bigvee_{m=1}^\infty \bigvee_{j=1}^{\nu^m} G(L_0)\wedge S^{m+1}\right),
$$
where $G(L_0)$ has discrete topology and $X\ltimes Y=(X\times Y)/(X\times\ast)$. Moreover, let $\calC(L_0)$ be the set of conjugation classes of $G(L_0)$. Then the quotient group
$$
\pi_n\left( \bigvee_{j=1}^\nu \calC(L_0)\ltimes S^2\right)\cong\pi_n\left(\bigvee_{j=1}^\nu S^2\right)\oplus\pi_n\left(\bigvee_{m=1}^\infty \bigvee_{j=1}^{\nu^m} \calC(L_0)\wedge S^{m+1}\right)
$$ gives the link invariants.
\end{enumerate}
\end{thm}
In assertion (2), the link invariant means the following: Given any element $\alpha \in \calA(L,L)\cong \pi_n(S^3)$, let $\tilde \alpha$ be an element in the intersection subgroup $A(L,L)$ that projects to $\alpha$ and let $l$ be a knot in $S^3\smallsetminus L$ that represents the element $\tilde \alpha$. Then the element $\alpha$ only depends on the ambient isotopy class of the $(n+1)$-link $L\cup l$. The meaning of link invariants in assertion (3) is similar.

Some consequences of (2) are as follows. Let $L$ be a strongly nonsplittable $n$-link. Let $\mathcal{L}^A$ be the set of ambient homotopy classes of the links $L\cup l$ with $l$ a knot in $S^3\smallsetminus L$ representing elements in the intersection subgroup $A(L,L)$. From assertion (2), there is a decomposition
$$
\mathcal{L}^A=\coprod_{\alpha\in \pi_n(S^3)} \mathcal{L}^A_{\alpha},
$$
where $\mathcal{L}^A_{\alpha}$ is given by those $L\cup l$ with $l$ representing $\alpha$. The connected sum operation on links induces an operation
$
\mathcal{L}^A_{\alpha}\times \mathcal{L}^A_{\beta}\longrightarrow \mathcal{L}^A_{\alpha+\beta}
$
for $\alpha,\beta\in \pi_n(S^3)$. Note the intersection subgroup $A(L,L)$ depends on the link group $G(L)$, but its quotient group $\calA(L,L)$ is independent on the choice of strongly nonsplittable $n$-links $L$. For the case $n=3$, the group $\calA(L,L)\cong \pi_3(S^3)=\Z$. For $n>3$, the group $\calA(L,L)\cong \pi_n(S^3)$ is a finite abelian group from the well-known results of Serre~\cite{Serre}. For instance from the table of homotopy groups in~\cite{Toda}, $\pi_4(S^3)=\pi_5(S^3)=\Z/2$ and $\pi_6(S^3)=\Z/12$. Up to the range that the homotopy groups are known (for instance the number of the link components of $L$ $\leq 64$), the groups $\calA(L,L)$ are fully determined by the above theorem. It should be pointed out that the homotopy groups $\pi_n(S^3)$ remains unknown for general $n$.

Note that the link group $G(L_0)$ always has (countably) infinite elements when $L_0\not=\emptyset$. From assertion (3), $\calA(L,L\smallsetminus L_0)$ contains the homotopy group $\pi_n(S^m)$ as summands with countably infinite occurrences for each $2\leq m\leq n$. In addition, the stable homotopy groups of spheres $\pi_n(S^m),$ $2\leq \frac{n+1}{2}\leq m\leq n$, occur  (countably) infinite times as summands in $\calA(L,L\smallsetminus L_0)$. The torsion free summands of the homotopy groups of the wedges of spheres can be understood through the Serre Theorem on the homotopy groups of spheres and the Hilton-Milnor Theorem on the decompositions of the loop space of the wedges of spaces. So the torsion free summands of $\calA(L,L\smallsetminus L_0)$ can be understood in this sense, while the torsion component of $\calA(L,L\smallsetminus L_0)$, which consists of those elements in $A(L,L\smallsetminus L_0)$ with certain powers lying in the symmetric commutator subgroup $A_S[L,L\smallsetminus L_0]$, remains as mystery.

Our method for obtaining Theorem~\ref{theorem1.1} is to consider the homotopy theory on the cubical diagrams of spaces induced by link complements through deleting link components. The classical Papakyriakopoulos Theorem~\cite{Pa} determines the homotopy type of link complements. For strongly nonsplittable pair of links $(L,L_0)$, the cubical diagram of the link complements obtained by removing the link components of $L\smallsetminus L_0$ has the special property that the spaces in the cubical diagram are all $K(\pi,1)$-spaces except the terminal space. Then the homotopy theory is able to establish the connections between the fundamental groups and the higher homotopy groups of the terminal space in these special cases. Technically the strongly nonsplittable hypothesis is necessary for applying the current results in homotopy theory on $n$-cubes to this problem. The determination of $\calA(L,L\smallsetminus L_0)$ remains open in general. In the case that $L\smallsetminus L_0$ is a $2$-link,  we are able to determine the group $\calA(L,L\smallsetminus L_0)$ (Theorem~\ref{theorem4.4}). For the case $L\smallsetminus L_0$ a $3$-link, certain canonical quotient of $\calA(L,L\smallsetminus L_0)$ can be determined in terms of the third homotopy groups (Theorem~\ref{theorem4.9}).

It is a natural question for comparing the homotopy-groups invariant $\calA(L,L_0)$ with other type of link invariants. We consider the Hopf link $L_n$ given as the pre-image of the $n$ distinct points in $S^2$ through the Hopf map $S^3\to S^2$. This gives a canonical example of strongly nonsplittable $n$-link. Theorem~\ref{theorem5.1} states that, for $n\geq4$, the group $A(L_n,L_n)$ has the trivial image in the homotopy link group in Milnor's sense and so $\pi_n(S^3)$-invariants can not be detected by Milnor's invariants for $n\geq4$.~\cite[Problem 1.96]{Kirby} asks to describe the links with vanishing Milnor's invariants.

We should point out that removing-component is a canonical operation on links or link diagrams and so it is possible to have other methods to study the group $\calA(L,L\smallsetminus L_0)$. We only investigate this object from the homotopy-theoretic views in this article. For highlighting the connections between links and the homotopy groups, we also provide a method how to construct homotopy group elements from $\calA(L,L_0)$ (Theorem~\ref{theorem6.5}). Our construction from the elements in $A(L, L\smallsetminus L_0)$ to $\pi_n(S^3\smallsetminus L_0)$ is the geometric analogue of the Massey products~\cite{Massey}, which may be called \textit{geometric Massey products}. Homologically there have been important applications of the Massey products to manifolds~\cite{DGMS, LS}. A direct consequence of Theorems~\ref{theorem1.1} and~\ref{theorem6.5} is that all homotopy groups of any higher dimensional spheres can be obtained from geometric Massey products on links.

The article is organized as follows. In section~\ref{section2}, we give a brief review on the homotopy type of link complements. The proof of Theorem~\ref{theorem1.1} is given in section~\ref{section3}. In section~\ref{section4}, we investigate the group $\calA(L,L')$ for $2$ or $3$-sublinks $L'$. We compare the homotopy-group invariants and Milnor's invariants in section~\ref{section5}. In section~\ref{section6}, we give some examples of links labeled by homotopy group elements as well as some remarks for constructing homotopy group elements from links.

\section{The Homotopy Type of Link Complements}\label{section2}
Let $L$ be a link in $S^3$. Suppose that $L$ is splittable. Then there is an embedding of $S^2$ into $S^3$ such that
$$
L\cong L'\sqcup L''
$$
with $L'$ and $L''$ located in different sides of $S^3\smallsetminus S^2$. It follows that there is a connected sum decomposition
$$
S^3\smallsetminus L \cong (S^3\smallsetminus L')\# (S^3\smallsetminus L'').
$$
By continuing this procedure, one gets a connected sum decomposition
\begin{equation}\label{equation2.1}
S^3\smallsetminus L \cong (S^3\smallsetminus L^{[1]})\#\cdots\# (S^3\smallsetminus L^{[\nu+1]})
\end{equation}
such that each $L^{[i]}$ is a nonsplittable sublink of $L$. The decomposition
\begin{equation}\label{equation2.2}
L\approx L^{[1]}\sqcup L^{[2]}\sqcup\cdots \sqcup L^{[\nu+1]}
\end{equation}
is called a \textit{complete splitting decomposition} of $L$. The number $\nu$ is called the \textit{splitting genus} of $L$, denoted by $\nu(L)$. Note that $L$ is nonsplittable if and only if $\nu(L)=0$.

The following classical theorem is due to Papakyriakopoulos~\cite[Theorem 1]{Pa}.
\begin{thm}[Papakyriakopoulos Theorem]\label{theorem2.1}
A link $L$ in $S^3$ is nonsplittable if and only if the link complement $S^3\smallsetminus L$ is a $K(\pi,1)$-space.
\end{thm}
Let $U$ be a nonempty proper open connected subset of the $3$-sphere $S^3$. \cite[Theorem 1]{Pa} states that $U$ is aspherical if and only if $S^3\smallsetminus U$ is nonsplittable, where the asphericity of $U$ is defined as $\pi_2(U)=0$. We give an elementary proof that $U$ is a $K(\pi,1)$-space if and only if $\pi_2(U)=0$ for any connected noncompact triangulated $3$-manifold $U$.
\begin{proof}
If $U$ is a $K(\pi,1)$, then $\pi_2(U)=0$ by definition of $K(\pi,1)$-spaces. Assume that $\pi_2(U)=0$. Let $\tilde U$ be the universal cover of $U$ with the covering map $q\colon \tilde U\to U$. Then $\tilde U$ is a connected noncompact $3$-manifold because otherwise $U=q(\tilde U)$ is compact which contradicts to that $U$ is open in $S^3$. Thus the third integral homology $H_3(\tilde U)=0$. Since $U$ is triangulated, so is $\tilde U$. From the fact that $$\pi_2(\tilde U)\cong \pi_2(U)=0$$ together with the fact that $\pi_1(\tilde U)=0$, the $3$-dimensional complex $\tilde U$ is homotopy equivalent to a wedge of $3$-spheres. It follows that $\tilde U$ is contractible because $H_3(\tilde U)=0$ and hence $U$ is a $K(\pi,1)$-space.
\end{proof}

\begin{prop}\label{proposition2.2}
Let $L$ be a link in $S^3$ with a complete splitting decomposition given in (\ref{equation2.2}). Then there is a homotopy decomposition
$$
S^3\smallsetminus L\simeq \left(\bigvee_{i=1}^{\nu}S^2 \right)\vee \left(\bigvee_{i=1}^{\nu+1} S^3\smallsetminus L^{[i]}\right).
$$
\end{prop}
\begin{proof}
Take an tubular neighborhood $V(L^{[i]})$ in $S^3$. Then  $S^3\smallsetminus V(L^{[i]})$ is a smooth compact manifold with nonempty boundary and $S^3\smallsetminus L^{[i]}\simeq S^3\smallsetminus V(L^{[i]})$. Given distinct points $q_1,\ldots,q_t$ in the interior of $S^3\smallsetminus V(L^{[i]})$, then
$$
S^3\smallsetminus (V(L^{[i]})\cup \{q_1,q_2,\ldots,q_t\})\simeq (S^3\smallsetminus (V(L^{[i]}))\vee \bigvee_{j=1}^t S^2_j,
$$
where $S^2_j$ is the boundary of a small ball around $q_i$. Now the assertion follows straightforward from the connected sum construction.
\end{proof}

\section{Proof of Theorem~\ref{theorem1.1}}\label{section3}
\subsection{$K(\pi,1)$-Partitions of Spaces}
Let $(X,X_0)$ be a pair of spaces. An  \textit{$n$-partition} of $X$ relative to $X_0$, denoted by $\mathbb{X}=(X;X_1,\ldots,X_n;X_0)$, means a sequence of subspaces $(X_1,\ldots,X_n)$ of $X$ such that
\begin{enumerate}
\item $X_0=X_i\cap X_j$ for each $1\leq i<j\leq n$ and
\item $X=\bigcup\limits_{i=1}^n X_i$.
\end{enumerate}
For any subset $I\subseteq \{1,\ldots,n\}$, let
$$
X_I=\bigcup_{i\in I}X_i,
$$
where $X_{\emptyset}=X_0$. A \textit{morphism}
$$
f\colon \mathbb{X}=(X;X_1,X_2,\ldots,X_n;X_0)\longrightarrow \mathbb{Y}=(Y;Y_1,Y_2,\ldots,Y_n;Y_0)
$$
between $n$-partitions of spaces means a (pointed continuous) map $f\colon X\to Y$ such that $f(X_i)\subseteq Y_i$ for $0\leq i\leq n$, where the basepoints for $X_i$ and $Y_i$ are chosen in $X_0$ and $Y_0$, respectively. An $n$-partition $(X;X_1,\ldots,X_n;X_0)$ is called \textit{cofibrant} if the inclusions
$$
X_I \hookrightarrow X_J
$$
are cofibrations for any $I\subseteq J\subseteq\{1,2,\ldots,n\}$. Note that for a cofibrant partition each union $X_I=\bigcup_{i\in I}X_i$ is the homotopy colimit of the diagram given by the inclusions $X_{I'}\hookrightarrow X_{I''}$ for $\emptyset\subseteq I'\subseteq I''\subsetneq I$, where $X_{\emptyset}=X_0$.

An $n$-partition $(X;X_1,\ldots,X_n;X_0)$ of $X$ is called an \textit{$K(\pi,1)$ $n$-partition} if the space $X_I$ is a path-connected $K(\pi,1)$-space for any proper subset $I\subsetneq \{1,2,\ldots,n\}$. Namely all subspaces $X_I$ except the total space $X=\bigcup_{i=1}^n X_i$ are $K(\pi,1)$-spaces. The spaces with $K(\pi,1)$-partitions have the following important property:

\begin{thm}\label{theorem3.1}
Let $(X;X_1,X_2,\ldots,X_n;X_0)$ be a cofibrant $K(\pi,1)$ $n$-partition with $n\geq 2$. Suppose that the inclusion $X_0 \to X_i$ induces an epimorphism on the fundamental groups for each $1\leq i\leq n$. Let $R_i$ be the kernel of $\pi_1(X_0)\to \pi_1(X_i)$ for $1\leq i\leq n$. Then
\begin{enumerate}
\item[(i)] For any proper subset $I=\{i_1,\ldots,i_k\}\subsetneq \{1,2,\ldots,n\}$,
$$
R_{i_1}\cap\cdots\cap R_{i_k}=[[R_{i_1}, R_{i_2}],\ldots,R_{i_k}]_S.
$$
\item[(ii)] For any $1<k\leq n$ and any subset $I=\{i_1,\ldots,i_k\}\subseteq \{1,2,\ldots,n\}$, there is an isomorphism of $\Z[\pi_1(X_0)]$-modules
$$
\rho_{\mathbb{X}}\colon \left.\left(\bigcap_{s=1}^k\left(R_{i_s}\cdot \prod_{j\in J}R_j\right)\right)\right/\left(\left([[R_{i_1}, R_{i_2}],\ldots,R_{i_k}]_S\right)\cdot \prod_{j\in J}R_j \right)\rTo^{\cong} \pi_k(X),
$$
where $J= \{1,2,\ldots,n\}-I$, $\pi_1(X_0)$ acts on each $R_i$ by conjugation and on $\pi_k(X)$ via the homomorphism $\pi_1(X_0)\to \pi_1(X)$ induced by the inclusion $X_0\hookrightarrow X$. Moreover the isomorphism is natural with respect to the morphisms of $n$-partitions.
In particular, there is a natural isomorphism of $\Z[\pi_1(X_0)]$-modules
$$
\rho_{\mathbb{X}}\colon (R_1\cap R_2\cap \cdots\cap R_n)/[[R_1,R_2],\ldots,R_n]_S\rTo^{\cong} \pi_n(X)$$
with respect to the morphisms of $n$-partitions.
\end{enumerate}
\end{thm}
Assertion (i) is given in~\cite[Theorem 1.3]{Li-Wu}. For assertion (ii),  was proved in ~\cite[Theorem 1.3]{Li-Wu} that
$$
 \left.\left(\bigcap_{s=1}^k\left(R_{i_s}\cdot \prod_{j\in J}R_j\right)\right)\right/\left(\left([[R_{i_1}, R_{i_2}],\ldots,R_{i_k}]_S\right)\cdot \prod_{j\in J}R_j \right)\cong \pi_k(X),
$$
as groups. The naturality is given in the following sense: Let $$\mathbb{X}=(X;X_1,X_2,\ldots,X_n;X_0)\textrm{ and }\mathbb{Y}=(Y;Y_1,Y_2,\ldots,Y_n;Y_0)$$ be two $n$-partitions satisfying the hypothesis in the theorem. Let $f\colon \mathbb{X}\to \mathbb{Y}$ be a morphism of $n$-partitions. Let $$R_i=\Ker(\pi_1(X_0)\to \pi_1(X_i))\textrm{ and }R'_i=\Ker(\pi_1(Y_0)\to \pi_1(Y_i))$$ for $1\leq i\leq n$. Then there is a commutative diagram
\begin{diagram}
 \left.\left(\bigcap_{s=1}^k\left(R_{i_s}\cdot \prod_{j\in J}R_j\right)\right)\right/\left(\left([[R_{i_1}, R_{i_2}],\ldots,R_{i_k}]_S\right)\cdot \prod_{j\in J}R_j \right)&\rTo^{\rho_{\mathbb{X}}}_{\cong} &\pi_k(X)\\
\dTo>{f_*}&&\dTo>{f_*}\\
 \left.\left(\bigcap_{s=1}^k\left(R'_{i_s}\cdot \prod_{j\in J}R'_j\right)\right)\right/\left(\left([[R'_{i_1}, R'_{i_2}],\ldots,R'_{i_k}]_S\right)\cdot \prod_{j\in J}R'_j \right)&\rTo^{\rho_{\mathbb{Y}}}_{\cong} &\pi_k(Y).
\end{diagram}
Here the homomorphism on the right column $f_*\colon \pi_k(X)\to \pi_k(Y)$ is induced by the map $f$ and the homomorphism on the left column is described as follows. The map $f\colon (X_i,X_0)\to (Y_i,Y_0)$ induces a commutative diagram
\begin{diagram}
R_i&\rInto & \pi_1(X_0)&\rOnto&\pi_1(X_i)\\
\dTo>{f_*}&&\dTo>{f_*}&&\dTo>{f_*}\\
R'_i&\rInto & \pi_1(Y_0)&\rOnto&\pi_1(Y_i)\\
\end{diagram}
and so it induces a group homomorphism $f_*$ in the left column of the above commutative diagram.

In order to see the naturality, the notion of $n$-cube of spaces is used. Let $\{0,1\}$ be the category with objects $\{0,1\}$ and the (non-identity) morphism given by the order $0<1$. Its $n$-fold Cartesian product is denoted by $\{0,1\}^n$. Similarly we have the defined category $\{-1,0,1\}^n$. Let $\mathrm{Top_*}$ denote the category of pointed spaces. A \textit{$n$-cube of spaces} $\mathbf{X}$ is a functor $\mathbf{X}\colon \{0,1\}^n\longrightarrow \mathrm{Top_*}$.
In other words, $\mathbf{X}$ is given by a $n$-cubical diagram of pointed spaces labeled by
$$\mathbf{X}(\epsilon)=X_{\epsilon_1,\ldots,\epsilon_n}$$
for $\epsilon=(\epsilon_1,\ldots,\epsilon_n)$ with $\epsilon_i=0,1$ and pointed maps
$$
X_{\epsilon_1,\ldots,\epsilon_n}\longrightarrow X_{\eta_1,\ldots,\eta_n}
$$
for $(\epsilon_1,\ldots, \epsilon_n)<(\eta_1,\ldots,\eta_n)$ in the category $\{0,1\}^n$. For an $n$-partition $\mathbb{X}=(X;X_1,\ldots,X_n;X_0)$, the associated $n$-cube of spaces $\mathbf{X}$ is given in the canonical way that $X_{(0,\ldots,0)}=X_0$, $X_{(0,\ldots,0,\stackrel{i}{1},0,\ldots,0)}=X_i$ and
$$
X_{\epsilon_1,\ldots,\epsilon_n}=\bigcup_{\epsilon_j=1}X_j
$$
with the maps in the cubical diagram given by the inclusions. According to~\cite[Chapter 3]{EH}, each $n$-cube of spaces $\mathbf{X}$ admits a natural embedding $\mathbf{X}\to \bar{\mathbf{X}}$ such that (i) each map $\mathbf{X}(\alpha)\to \bar{\mathbf{X}}(\alpha)$ is a homotopy equivalence with natural homotopy inverse and (ii) $\bar{\mathbf{X}}$ is \textit{fibrant} in the sense that for each $\epsilon\in \{0,1\}^n$ the canonical map $\bar {\mathbf{X}}(\epsilon)\to \lim_{\alpha>\epsilon}\bar {\mathbf{X}}(\alpha)$ is a fibration. Moveover from~\cite{Steiner} such a fibrant $n$-cube $\bar{\mathbf{X}}$ may be extended to an $n$-cube of fibrations in the sense of a functor from the category $\{-1,0,1\}^n$ to pointed spaces, also denoted as $\bar{\mathbf{X}}$, such that, for each $1\leq k\leq n$ and and $\epsilon\in \{-1,0,1\}^n$, $\bar{\mathbf{X}}(\epsilon_1,\ldots,\epsilon_{k-1},-1,\epsilon_{k+1},\ldots,\epsilon_n)$ is the fibre of the fibration
$$
\bar{\mathbf{X}}(\epsilon_1,\ldots,\epsilon_{k-1},0,\epsilon_{k+1},\ldots,\epsilon_n)\longrightarrow \bar{\mathbf{X}}(\epsilon_1,\ldots,\epsilon_{k-1},1,\epsilon_{k+1},\ldots,\epsilon_n).
$$
The construction $\mathbf{X}\mapsto \bar{\mathbf{X}}$ gives a functor from the category of $n$-cubes of spaces to the $n$-cubes of fibrations.

\begin{proof}[Proof of Theorem~\ref{theorem3.1}]
We only need to show that the isomorphism in assertion (2) is a natural isomorphism of $\Z[\pi_1(X_0)]$-modules.
Let $\mathbf{X}$ be the $n$-cube of spaces induced by the $n$-partition $\mathbb{X}=(X;X_1,\ldots,X_n;X_0)$.
Observe that~\cite[Theorem 1.3]{Li-Wu} is obtained by showing that the connectivity hypothesis in~\cite[Theorem 1]{EM} holds for the subgroups $R_1,R_2,\ldots,R_n$. On the other hand, ~\cite[Theorem 1]{EM} is obtained by inspecting the homotopy exact sequences of the fibrations in the $n$-cube of the fibrations $\bar{\mathbf{X}}$. From the functor $\mathbb{X}\mapsto \bar{\mathbf{X}}$, the isomorphism in (2) is a natural isomorphism of $\Z[\pi_1(X_0)]$-modules and hence the result.
\end{proof}

\subsection{Proof of Theorem~\ref{theorem1.1}}
Let $M=S^3\smallsetminus L_0$ and let $L\smallsetminus L_0=\{l_1,\ldots,l_n\}$. Let
$$
M_i=M\smallsetminus \{l_1,\ldots,l_{i-1},l_{i+1},\ldots,l_n\}
$$
for $1\leq i\leq n$ with
$$
M_0=M\smallsetminus (L\smallsetminus L_0)=S^3\smallsetminus L.
$$
Then $\mathbb{M}=(M;M_1,\ldots, M_n;M_0)$ is a cofibrant $n$-partition of $M$ relative to $M_0$. We need to check the $K(\pi,1)$-hypothesis. Let $I\subsetneq \{1,2,\ldots,n\}$ be a proper subset. Then
$$
K=\{l_j \ | \ j\not\in I\}
$$
is a nonempty sublink of $L\smallsetminus L_0$ and so
$$
M_I=\bigcup_{i\in I} M_i= M\smallsetminus K=S^3\smallsetminus(L_0\cup K).
$$
is $K(\pi,1)$ space by Theorem~\ref{theorem2.1} because $L_0\cup K$ ($\supsetneq L_0)$ is nonsplittable by the definition of strong nonsplittablity. Clearly
$$
\pi_1(M_0)\longrightarrow \pi_1(M_i)
$$
is an epimorphism because one can deform any loop in $M_i$ such that it does not intersect with $l_i$. Thus we can apply Theorem~\ref{theorem3.1}. Note that
$$
A(L,l_i)=\Ker(\pi_1(M_0)\to \pi_1(M_i))
$$
is the group $R_i$ in Theorem~\ref{theorem3.1}.

(1). If $L'\subsetneq L\smallsetminus L_0$, the assertion that $\calA(L,L')=0$ follows from Theorem~\ref{theorem3.1}(i). If $L'=L\smallsetminus L_0$ and $L_0$ is nonempty and nonsplittable, then $M$ is a $K(\pi,1)$-space. the assertion then follows from Theorem~\ref{theorem3.1}(ii). The first part of (2) is also a direct consequence of Theorem~\ref{theorem3.1}(ii).

For the first part of (3), from Theorem~\ref{theorem3.1}(ii),
\begin{equation}\label{equation3.1}
\calA(L,L\smallsetminus L_0)\cong \pi_n(M).
\end{equation}
By Proposition~\ref{proposition2.2},
$$
M\simeq \left(\bigvee\limits_{i=1}^{\nu}S^2 \right)\vee \left(\bigvee\limits_{i=1}^{\nu+1} S^3\smallsetminus L^{[i]}\right),
$$
where $L_0=\coprod_{i=1}^{\nu+1}L^{[i]}$ is the complete splitting decomposition of $L_0$. Since each $S^3\smallsetminus L^{[i]}$ is a $K(\pi,1)$-space, so is the wedge sum $\bigvee_{i=1}^{\nu+1} S^3$. Thus
\begin{equation}\label{equation3.2}
M\simeq K(G,1)\vee \bigvee_{i=1}^{\nu}S^2,
\end{equation}
where $G=\pi_1(M)=G(L_0)$. By taking the universal covering, we have
\begin{equation}\label{equation3.3}
\tilde M\simeq G(L_0)\ltimes \left(\bigvee_{i=1}^{\nu}S^2\right)=\bigvee_{i=1}^\nu G(L_0)\ltimes S^2.
\end{equation}
It follows that
$$
\pi_n(M)\cong \pi_n\left(\bigvee_{i=1}^\nu G(L_0)\ltimes S^2\right),
$$
which gives the first isomorphism in assertion (3). Let $X=\bigvee_{i=1}^\nu S^1$. From the Hilton-Milnor Theorem~\cite{Gray, Hilton, Milnor2} (also see Whitehead's book~\cite[Sections 6-7, Chapter XI]{Whitehead}) together with the suspension splitting theorem of loop suspensions~\cite[Corollary 2.11, p.335]{Whitehead},
$$
\begin{array}{rcl}
\Omega\left( \bigvee\limits_{i=1}^{\nu}G(L_0)\ltimes S^2\right)&=& \Omega\Sigma\left(X\vee  (G(L_0)\wedge X)\right)\\
&\simeq& \Omega\Sigma X  \times \Omega \left( \Sigma (G(L_0)\wedge X)\vee \Sigma (G(L_0)\wedge X)\wedge \Omega \Sigma X\right)\\
&=& \Omega \Sigma X\times \Omega\left(G(L_0)\wedge (\Sigma X \vee (\Sigma X\wedge\Omega \Sigma X))\right)\\
&\simeq&\Omega \Sigma X \times \Omega \left(G(L_0)\wedge \left(\Sigma X\vee \Sigma X\wedge\left(\bigvee\limits_{k=1}^\infty X^{\wedge k}\right)\right)\right)\\
&=&\Omega \Sigma X \times \Omega \Sigma \left(G(L_0)\wedge \bigvee_{m=1}^\infty X^{\wedge m}\right)\\
&=&\Omega \left(\bigvee\limits_{i=1}^{\nu}S^2\right) \times \Omega \left(\bigvee\limits_{m=1}^\infty \bigvee\limits_{i=1}^{\nu^m} G(L_0)\wedge S^{m+1}\right).\\
\end{array}
$$
It follows that
$$
\pi_n\left(\bigvee_{i=1}^\nu G(L_0)\ltimes S^2\right)\cong \pi_n\left(\bigvee\limits_{i=1}^{\nu}S^2\right)\bigoplus \pi_n\left(\bigvee\limits_{i=1}^{\nu}S^2\right) \times \Omega \left(\bigvee\limits_{m=1}^\infty \bigvee\limits_{i=1}^{\nu^m} G(L_0)\wedge S^{m+1}\right),
$$
which gives the second isomorphism. Similarly
$$
\pi_n\left( \bigvee_{j=1}^\nu \calC(L_0)\ltimes S^2\right)\cong\pi_n\left(\bigvee_{j=1}^\nu S^2\right)\oplus\pi_n\left(\bigvee_{m=1}^\infty \bigvee_{j=1}^{\nu^m} \calC(L_0)\wedge S^{m+1}\right).
$$

Now we prove the link invariant property in (2) and (3). Let $l$ be a knot that representing an element $\alpha\in A(L,L\smallsetminus L_0)$ in following sense: Choose a path $\lambda$ in $S^3\smallsetminus L$ from the basepoint to a point in $l$. Then the loop homotopy class of the path product $\lambda\ast l\ast \lambda^{-1}$ give the element $\alpha\in A(L,L\smallsetminus L_0)$. Since $A(L,L_0)$ is a normal subgroup of $\pi_1(S^3\smallsetminus L)$, a different choice of the paths $\lambda$ gives a conjugation of the element $\alpha$ and so the knot $l$ determines a unique conjugation class of $\alpha$.

Assume that $\bar L\cup \bar l$ is a link in $S^3$ which is ambient isotopic to $L\cup l$. There exists an ambient isotopy $h_t\colon S^3\to S^3$ with $h_0=\id_{S^3}$ and $h_1(L\cup l)=\bar L\cup \bar l$. The ambient isotopy $h_t$ induces a homeomorphism given by
\begin{equation}\label{equation3.4}
H\colon S^3\times I\longrightarrow S^3\times I\quad H(x,t)=(h_t(x),t),
\end{equation}
where $I=[0,1]$. Let $X_0=(S^3\times I)\smallsetminus H(L\times I)$ and
$$
X_i=(S^3\times I)\smallsetminus H((L\smallsetminus {l_i})\times I)
$$
with $X=(S^3\times I)\smallsetminus H(L_0\times I)$. The homeomorphism in~(\ref{equation3.4}) induces a homeomorphism
\begin{equation}\label{equation3.5}
H\colon M_J\times I\longrightarrow X_J=\bigcup_{j\in J} X_j
\end{equation}
for any subset $J\subseteq \{1,2,\ldots,n\}$. Thus $\mathbb{X}=(X;X_1,\ldots,X_n;X_0)$ is a cofibrant $K(\pi,1)$ $n$-partition of $X$ relative to $X_0$. Let $(R_i;a)$ be the kernel of
$$
\pi_1(X_0;a)\longrightarrow \pi_1(X_i;a)
$$
for a choice of the basepoint $a\in X_0$.

Let $\bar M_0=S^3\smallsetminus \bar L$ and
$$
\bar M_i=S^3\smallsetminus \{\bar l_1,\ldots,\bar l_{i-1},\bar l_{i+1},\ldots,\bar l_n\}
$$
with $\bar M=S^3\smallsetminus \bar L_0$. From the homeomorphism $h_1\colon M_J\to \bar M_J$, $$\bar{\mathbb{M}}=(\bar M; \bar M_1,\ldots,\bar M_n;\bar M_0)$$ is a cofibrant $K(\pi,1)$ $n$-partition of $\bar M$ relative to $\bar M_0$ and the knot $\bar l$ in $\bar M_0$ represents a conjugation class of an element $\bar \alpha\in A(\bar L, \bar L\smallsetminus \bar L_0)$. Let
$$
f\colon M_J=M_J\times \{0\}\hookrightarrow X_J\textrm{ and } g\colon \bar M_J=\bar M_J\times \{1\}\hookrightarrow X_J
$$
be the canonical inclusion. We identify $M_J$ with $f(M_J)$ and  $\bar M_J$ with $g(M_J)$.

By homeomorphism~(\ref{equation3.5}), both $M_J$ and $\bar M_J$ are strong deformation retracts of $X_J$. By choosing a basepoint $x_1\in \bar M_0\times \{1\}$, the inclusion $g\colon \bar M_J\to X_J$ induces a canonical isomorphism
$$
g_*\colon A(\bar L,\bar L\smallsetminus \bar L_0)\cong (R_1; x_1)\cap (R_2;x_1)\cap\cdots\cap (R_n;x_1).
$$
Similarly there is a canonical isomorphism
$$
f_*\colon A(L,L\smallsetminus L_0)\cong (R_1; x_0)\cap (R_2;x_0)\cap\cdots\cap (R_n;x_0)
$$
for a basepoint $x_0\in M_0\times\{0\}$.
Let $\zeta$ be any path in $X_0$ from $x_0$ to $x_1$ and let
$$
\zeta_*\colon \pi_1(X_0;x_0)\rTo^{\cong} \pi_1(X_0;x_1) \textrm{ and } \zeta_*\colon \pi_n(X;x_0)\rTo^{\cong} \pi_n(X;x_1)
$$
be the induced isomorphisms of the path $\zeta$. Note that the resulting isomorphism
$$
\zeta_*\colon \pi_n(X;x_0)\otimes_{\Z[\pi_1(X_0;x_0)]}\Z\rTo^{\cong} \pi_n(X;x_1)\otimes_{\Z[\pi_1(X_0;x_1)]}\Z
$$
is independent on the choices of the paths $\zeta$ from $x_0$ to $x_1$. From deformation $h_t(l)\subseteq X_0$ with $0\leq t\leq 1$, $\zeta(f_*(\alpha))$ is conjugate to $g_*(\bar\alpha)$ in $\pi_1(X_0;x_1)$. Thus
$$
\zeta_*(\rho_{\mathbb{X}}(f_*(\alpha)))\equiv \rho_{\mathbb{X}}(g_*(\bar\alpha))
$$
in the quotient group $\pi_n(X;x_1)\otimes_{\Z[\pi_1(X_0;x_1)]}\Z$. By the naturality in Theorem~\ref{theorem3.1}(ii), we have
$$
\rho_{\mathbb{X}}(f_*(\alpha))=f_*(\rho_{\mathbb{M}}(\alpha))\textrm{ and } \rho_{\mathbb{X}}(g_*(\bar\alpha))=g_*(\rho_{\bar{\mathbb{M}}}(\bar\alpha)).
$$
It follows that
\begin{equation}\label{equation3.6}
\zeta_*(f_*(\rho_{\mathbb{M}}(\alpha)))\equiv g_*(\rho_{\bar{\mathbb{M}}}(\bar\alpha)).
\end{equation}
in the quotient group $\pi_n(X;x_1)\otimes_{\Z[\pi_1(X_0;x_1)]}\Z$.

\noindent\textbf{Case 1.} $\L_0=\emptyset$ as in assertion (2). Then
$$
X=(S^3\times I)\smallsetminus H(L_0\times I)=S^3\times I
$$
as $L_0=\emptyset$. The maps $f\colon M=S^3\to X=S^3\times I$ and $g\colon \bar M=S^3\to  X\times I$ are the inclusions of $S^3\times \{0\}$ and $S^3\times\{1\}$ into the cylinder $S^3\times I$, respectively. Note that $X=S^3\times I$ is simply connected. The action of $\pi_1(X_0)$ on $\pi_n(X)$ is trivial. Moreover $\pi_n(X)$ is independent on the choice of the basepoints. From equation~(\ref{equation3.6}), we have the same element given in the homotopy group under the canonical identification
$$
\pi_n(S^3)\rEq^{f_*} \pi_n(S^3\times I)\rEq^{g_*} \pi_n(S^3).
$$
Note that this element is independent on the choice of the ambient isotopy $h_t$ because the final maps $f\colon S^3\to S^3\times I$ and $g\colon S^3\to S^3\times I$ are independent on $H$. This proves the link invariant property in assertion (2).

\noindent\textbf{Case 2.} $L_0\not=\emptyset$ and splittable with splitting genus $\nu$ as in assertion (3). By Proposition~\ref{proposition2.2},
$$
X\simeq S^3\smallsetminus L_0\simeq K(G(L_0),1)\vee \bigvee_{j=1}^\nu S^2.
$$
Let $\pi\colon  \widetilde{K(G(L_0),1)}\to K(G(L_0),1)$ be the universal covering. Then the universal covering of $K(G(L_0),1)\vee \bigvee_{j=1}^\nu S^2$ is the union
$$
\widetilde{K(G(L_0),1)}\cup \left(G(L_0)\times\left(\bigvee_{j=1}^\nu S^2\right)\right)
$$
by identifying the points $\pi^{-1}(\ast)=G(L_0)$ with the corresponding points in the subspace $G(L_0)\times\ast$ of $G(L_0)\times\left(\bigvee_{j=1}^\nu S^2\right)$. Thus the universal covering $\tilde X$ is $G(L_0)$-equivariantly homotopy equivalent to the space
$$
G(L_0)\ltimes\left(\bigvee_{j=1}^\nu S^2\right)
$$
with $G(L_0)$-action is given by
$$
g\cdot (h,x)=(ghg^{-1},x)
$$
for $g\in G(L_0)$ and $(h,x)\in G(L_0)\ltimes\left(\bigvee_{j=1}^\nu S^2\right)$. It follows that $\pi_1(X)$ acts trivially on the quotient space
$$
\calC(L_0)\ltimes \left(\bigvee_{j=1}^\nu S^2\right).
$$
From equation ~(\ref{equation3.6}), we have
$$
f_*(\rho_{\mathbb{M}}(\alpha)))= g_*(\rho_{\bar{\mathbb{M}}}(\bar\alpha))
$$
in the homotopy group
$$
\pi_n\left(\calC(L_0)\ltimes \left(\bigvee_{j=1}^\nu S^2\right)\right)
$$
under the canonical identification induced by $f\colon S^3\smallsetminus L_0\to X$ and $g\colon S^3\smallsetminus \bar L_0\to X$.
To finish the proof, we have to show that the above identification between $\rho_{\mathbb{M}}(\alpha)$ and $\rho_{\bar{\mathbb{M}}}(\bar\alpha)$ is independent on the choices of the ambient isotopy $h_t$. Let $h'_t$ be another ambient isotopy between $L\cup l$ and $\bar L\cup \bar l$. Consider the homeomorphism
$$
H'\colon S^3\times I\longrightarrow S^3\times I\quad H'(x,t)=(h'_t(x),t).
$$
Then there is a commutative diagram
\begin{diagram}
S^3\smallsetminus \bar L_0&\rEq& S^3\smallsetminus \bar L_0\\
\simeq \dInto>{g}&&\simeq \dInto>{g}\\
(S^3\times I)\smallsetminus H(L_0\times I)&\rTo^{H'\circ H^{-1}}& (S^3\times I)\smallsetminus H'(L_0\times I)\\
\simeq \uInto>{f} &&\simeq \uInto>{f}\\
S^3\smallsetminus L_0&\rEq& S^3\smallsetminus L_0.\\
\end{diagram}
By applying $\pi_n(\, )$ to the above diagram, the identification $f_*(\rho_{\mathbb{M}}(\alpha))=g_*(\rho_{\bar{\mathbb{M}}}(\bar\alpha))$
in the homotopy group
$$
\pi_n\left(\calC(L_0)\ltimes \left(\bigvee_{j=1}^\nu S^2\right)\right)
$$
is independent on the choice of the ambient isotopy $h_t$. This finishes the proof of Theorem~\ref{theorem1.1}.\hfill $\Box$

\section{The Groups $\calA(L,L')$ for Sub $n$-Link $L'$ With $n\leq 3$.}\label{section4}
\subsection{Lemmas on Homotopy $n$-Pushouts}
Recall that an \textit{$n$-corner} is a functor $\mathbf{X}$ from the category $\{0,1\}^n\smallsetminus \{(1,\ldots,1)\}$ to (pointed) spaces. Namely $\mathbf{X}$ is a cubical homotopy commutative diagram of (pointed) spaces without counting the terminal space $X_{(1,\ldots,1)}$. In $n=1$, $\mathbf{X}=X_0$ is a single space. The diagrams of $2$-corners and $3$-corners are pictured as follows:
\begin{diagram}
\begin{diagram}
X_{0,0}&\rTo& X_{0,1}\\
\dTo&&\\
X_{1,0}&&\\
\end{diagram}
&
\begin{diagram}
X_{0,0,0}&     & \rTo   &     & X_{1,0,0}&    &\\
        &\rdTo&         &     &\dDashto      &\rdTo&\\
\dTo    &     &X_{0,1,0}&\rTo &          &       &X_{1,1,0}\\
        &     &     &     &     &       &\\
X_{0,0,1}& \rDashto   &     &     &X_{1,0,1}&&\\
         &\rdTo& \dTo    & &&&\\
         &     &X_{1,1,0}& &&&\\
\end{diagram}\\
\end{diagram}
An $n$-corner $\mathbf{X}$ is called \textit{cofibrant} if for every $$ \epsilon=(\epsilon_1,\ldots,\epsilon_n)\in \{0,1\}^n\smallsetminus\{(1,\ldots,1)\}$$ the canonical map
$$
\colim\{X_{\eta}\ | \ \eta<\epsilon \}\longrightarrow \colim  \colim\{X_{\eta} \ | \ \eta\leq \epsilon\}=X_{\epsilon}
$$
is a cofibration. The homotopy colimit of an $n$-corner $\mathbf{X}$ is called a \textit{homotopy $n$-pushout}, denoted by $\hocolim \mathbf{X}$. Let $\mathbf{X}_{\mathrm{top}}$ be the full sub diagram of $\mathbf{X}$ consisting of the spaces in the top $(n-1)$-cubical diagram except the terminal space, namely the object of $\mathbf{X}_{\mathrm{top}}$ consists of the spaces
$$
X_{\epsilon_1,\ldots,\epsilon_{n-1},0}
$$
with $(\epsilon_1,\ldots,\epsilon_{n-1},0)\not=(1,\ldots,1,0)$ and the maps induced from that in $\mathbf{X}$. Let $\mathbf{X}_{\mathrm{bottom}}$ be the full sub diagram of $\mathbf{X}$ consisting of the spaces
$$
X_{\epsilon_1,\ldots,\epsilon_{n-1},1}
$$
with the maps induced from that in $\mathbf{X}$. The vertical maps in $\mathbf{X}$ induces a canonical morphism of $(n-1)$-corners
\begin{equation}\label{equation4.1}
\mathbf{f}\colon \mathbf{X}_{\mathrm{top}}\longrightarrow \mathbf{X}_{\mathrm{bottom}}
\end{equation}
and so a map
$$
\hocolim\mathbf{f}\colon \hocolim \mathbf{X}_{\mathrm{top}}\longrightarrow \hocolim\mathbf{X}_{\mathrm{bottom}}.
$$
By looking at the top $(n-1)$-cubical sub diagram of $\mathbf{X}$, there is a canonical map
\begin{equation}\label{equation4.2}
g\colon \hocolim\mathbf{X}_{\mathrm{top}}\longrightarrow X_{1,\ldots,1,0}
\end{equation}
because $X_{1,\ldots,1,0}$ is the terminal space in the top $(n-1)$-cubical diagram.

\begin{lem}\label{lemma4.1}
Let $X$ be an $n$-corner and let $f$ and $g$ be given as above. Then there is a homotopy pushout diagram
\begin{diagram}
\hocolim\mathbf{X}_{\mathrm{top}}&\rTo^g& X_{1,\ldots,1,0}\\
\dTo>{\hocolim\mathbf{f}}&&\dTo\\
\hocolim\mathbf{X}_{\mathrm{bottom}}&\rTo&\hocolim\mathbf{X}
\end{diagram}
\end{lem}
\begin{proof}
From the construction in~\cite{Vogt}, every $n$-corner is equivalent to a cofibrant $n$-corner. Thus we may assume that $\mathbf{X}$ is a cofibrant $n$-corner in which each map in the diagram is the inclusion of subspaces that are cofibrations. By~\cite[Proposition 1.6]{Goodwillie}, the canonical map $\hocolim \mathbf{Y}\to \colim \mathbf{Y}$ is a homotopy equivalence for any cofibrant $m$-corner $\mathbf{Y}$. Observe that
$$
\colim\mathbf{X}=\bigcup_{(\epsilon_1,\ldots,\epsilon_n)} X_{\epsilon_1,\ldots,\epsilon_n}.
$$
Since
$$
X_{\epsilon_1,\ldots,\epsilon_{n-1},0}\subseteq X_{\epsilon_1,\ldots,\epsilon_{n-1},1}
$$
for any $(\epsilon_1,\ldots,\epsilon_{n-1})\in \{0,1\}^n$ with $(\epsilon_1,\ldots,\epsilon_{n-1})\not=(1,\ldots,1)$,
$$
\colim\mathbf{X}=\left(\bigcup_{(\epsilon_1,\ldots,\epsilon_{n-1})\not=(1,\ldots,1)} X_{\epsilon_1,\ldots,\epsilon_{n-1},1}\right)\cup X_{1,\ldots,1,0}=\colim\mathbf{X}_{\mathrm{bottom}}\cup X_{1,\ldots,1,0}
$$
with the intersection
$$
\colim\mathbf{X}_{\mathrm{bottom}}\cap X_{1,\ldots,10}=\left(\bigcup_{(\epsilon_1,\ldots,\epsilon_{n-1})\not=(1,\ldots,1)} X_{\epsilon_1,\ldots,\epsilon_{n-1},0}\right)=\colim\mathbf{X}_{\mathrm{top}}
$$
and hence the result.
\end{proof}

Let $\mathbf{f}\colon \mathbf{X}\to \mathbf{Y}$ be a morphism of $n$-corners. Then we have the \textit{mapping cone} $\mathbf{C}_{\mathbf{f}}$ which is the $n$-corner with spaces given by the mapping cone of
the maps $f_{\epsilon}\colon X_{\epsilon}\to Y_{\epsilon}$ for $\epsilon\in \{0,1\}^n\smallsetminus\{(1,\ldots,1)\}$ with the canonical induced maps in the diagram.

\begin{lem}\label{lemma4.2}
Let $\mathbf{f}\colon \mathbf{X}\to \mathbf{Y}$ be a morphism of $n$-corners. Then there is a cofibre sequence
$$
\hocolim \mathbf{X}\longrightarrow \hocolim\mathbf{Y}\longrightarrow \hocolim\mathbf{C}_{\mathbf{f}}.
$$
\end{lem}
\begin{proof}
We assume that both $\mathbf{X}$ and $\mathbf{Y}$ are cofibrant $n$-corners. Moreover we may assume that the morphism $\mathbf{f}\colon \mathbf{X}\to \mathbf{Y}$ is a cofibration in the sense that
each map
$$
f_{\epsilon}\colon X_{\epsilon}\longrightarrow Y_{\epsilon}
$$
is a cofibration for each $\epsilon$ (for instance we can replace each $Y_{\epsilon}$ by the mapping cylinder of $f_{\epsilon}\colon X_{\epsilon}\to Y_{\epsilon}$). Then $\mathbf{C}_{\mathbf{f}}$ is equivalent to the cofibrant $n$-corner $\mathbf{Y}/\mathbf{X}=\{Y_{\epsilon}/X_{\epsilon}\}$. Thus
$$
\hocolim\mathbf{C}_{\mathbf{f}}\simeq \hocolim \mathbf{Y}/\mathbf{X}\simeq \colim \mathbf{Y}/\mathbf{X}=(\colim\mathbf{Y})/(\colim\mathbf{X})
$$
and hence the result.
\end{proof}

\subsection{The Groups $\calA(L,L')$ for Sub $2$-Links $L'$}
Let $L$ be any $n$-link in $S^3$ with splitting genus of $\nu\geq 0$. If $\nu\geq 1$, consider the connected sum decomposition~(\ref{equation2.1}). Let
$$
f_i\colon S^2\hookrightarrow S^3\smallsetminus L=(S^3\smallsetminus L^{[1]})\# (S^3\smallsetminus L^{[2]})\#\cdots\# (S^3\smallsetminus L^{[\nu+1]})
$$
be the inclusion of the $2$-sphere into the $i\,$th separating $2$-sphere in the decomposition for $1\leq i\leq \nu$. Let
$$
\widehat{S^3\smallsetminus L}=(S^3\smallsetminus L)\cup_{f_1}D^3\cup_{f_2}D^3\cup\cdots\cup_{f_\nu}D^3
$$
be obtained by attaching $\nu$ $3$-cells. The space $\widehat{S^3\smallsetminus L}$ is no longer a $3$-manifold but from Proposition~\ref{proposition2.2}, we have
$$
\widehat{S^3\smallsetminus L}\simeq K(G(L),1)
$$
which is a model for the classifying space of the link group $G(L)$.

Let $d_iL$ be the sublink of $L$ by removing the $i\,$th link component of $L$. For a subset $I=\{i_1,\ldots,i_k\}\subseteq \{1,\ldots,n\}$, let $d_IL$ be the sublink of $L$ by removing $i_1\,\ldots\,i_k\,$th link components of $L$. Then inclusions of $S^3\smallsetminus L$ into  $S^3\smallsetminus d_iL$ and $S^3\smallsetminus d_jL$ induces a $2$-corner
\begin{diagram}
\widehat{S^3\smallsetminus L}&\rTo& \widehat{S^3\smallsetminus d_iL}\\
\dTo&&\\
\widehat{S^3\smallsetminus d_jL}.&&\\
\end{diagram}
Let $X^L_{i,j}$ be the homotopy pushout of the above $2$-corner. The homotopy type of $X^L_{i,j}$ is as follows. A link component $l$ of $L$ is called \textit{splittable} in $L$ if $L=(L\smallsetminus l)\cup l$ is a splittable decomposition.

\begin{lem}\label{lemma4.3}
Let $L=\{l_1,\ldots,l_n\}$ be an $n$-link in $S^3$ with $n\geq 2$ and let $1\leq i\not=j\leq n$.
\begin{enumerate}
\item If $n=2$ and $\nu(L)=0$, then $X^L_{1,2}\simeq S^3$.
\item If $n=2$ and $\nu(L)=1$, then $X^L_{1,2}\simeq\ast$.
\item Let $n>2$.
\begin{enumerate}
\item[(i)] Suppose that $l_i$ and $l_j$ are linked together lying in the same nonsplittable components of $L$ with the property that there are no nonsplittable components of $d_{i,j}L$ that can be linked with both $l_i$ and $l_j$. Then
    $$X^L_{i,j}\simeq \widehat{S^3\smallsetminus d_{i,j}L}\vee S^3$$
with $\nu(d_{i,j}L)-\nu(d_iL)-\nu(d_jL)+\nu(L)=-1$.
\item[(ii)] Suppose that $l_i$ and $l_j$ are linked together lying in the same nonsplittable components of $L$ with the property that there are two or more nonsplittable components of $d_{i,j}L$ that can be linked with both $l_i$ and $l_j$. Then
$$
X^L_{i,j}\simeq \widehat{S^3\smallsetminus d_{i,j}L}\vee \bigvee_{k=1}^{\nu(d_{i,j}(L))-\nu(d_i(L))-\nu(d_j(L))+\nu(L)}S^2
$$
with $\nu(d_{i,j}L)-\nu(d_iL)-\nu(d_jL)+\nu(L)>0$.
\item[(iii)] In the remaining cases,
$$
X^L_{i,j}\simeq \widehat{S^3\smallsetminus d_{i,j}L}\simeq K(G(d_{i,j}L),1)
$$
is a $K(\pi,1)$-space with $\nu(d_{i,j}L)-\nu(d_iL)-\nu(d_jL)+\nu(L)=0$.
\end{enumerate}
\end{enumerate}
\end{lem}
An example of $3$-link for the case in assertion (3i) is to take a trivial $2$-link labeled by $l_i$ and $l_k$ with adding the third component $l_j$ linked with both $l_i$ and $l_k$. An example of $4$-link for the case in assertion (3ii) is to take a trivial $2$-link with adding two more components $l_i$ and $l_j$ linked with both components of the trivial $2$-link.
\begin{proof}
(1). In this case, the links $L$, $d_1L$ and $d_2L$ are not splittable. Thus $$\widehat{S^3\smallsetminus d_IL}\simeq S^3\smallsetminus d_IL$$ for $I=\emptyset, 1,2$ and so
$$
X^L_{1,2}\simeq (S^3\smallsetminus d_1L)\cup (S^3\smallsetminus d_2L)=S^3.
$$

(2). In this case, $S^3\smallsetminus L \simeq S^2\vee \widehat{S^3\smallsetminus L}$ and $S^3\smallsetminus d_iL \simeq \widehat{S^3\smallsetminus d_iL}$ for $i=1,2$. Note that $S^3$ is the homotopy pushout of $2$-corners $$\ast\longleftarrow S^2\longrightarrow \ast\textrm{ and  }S^3\smallsetminus d_1L\longleftarrow S^3\smallsetminus L\longrightarrow S^3\smallsetminus d_2L.$$ By Lemma~\ref{lemma4.2}, there is a cofibre sequence
$$
S^3\rTo^{\simeq} S^3\rTo X^L_{1,2}
$$
and so $X^L_{1,2}\simeq\ast$.

(3). The proof is given as case-by-case, where assertions (3i) and (3ii) follow from Sub-Sub-Case 3.2.1 and Sub-Sub-Case 3.2.3, respectively.

\noindent\textbf{Case 1.} \textit{Both $l_i$ and $l_j$ are splittable in $L$}. Then
$$
S^3\smallsetminus L\cong (S^3\smallsetminus d_{i,j}L)\# (S^3\smallsetminus l_i)\# (S^3\smallsetminus l_j)
$$
with
$$
\begin{array}{rcl}
\nu(L)&=&\nu(d_{i,j}L)+2,\\
\nu(d_iL)&=&\nu(d_{i,j}L)+1,\\
\nu(d_jL)&=&\nu(d_{i,j}L)+1.\\
\end{array}
$$
Note that $S^3\smallsetminus d_{i,j}L$ is the pushout of
\begin{equation}\label{equation4.3}
S^3\smallsetminus d_iL\longleftarrow S^3\smallsetminus L\longrightarrow S^3\smallsetminus d_jL.
\end{equation}
Let $A$ be the homotopy pushout of
$$
\bigvee_{k=1}^{\nu(d_iL)=\nu(d_{i,j}L)+1}S^2\longleftarrow \bigvee_{k=1}^{\nu(L)=\nu(d_{i,j}L)+2}S^2\longrightarrow \bigvee_{k=1}^{\nu(d_{j}L=\nu(d_{i,j}L)+1)} S^2
$$
by taking out the $2$-spheres in the above $2$-corner. Then
$$
A\simeq \bigvee_{k=1}^{\nu(d_{i,j}L)}S^2,
$$
which is homotopy equivalent to the union of the separating $2$-spheres in $S^3\smallsetminus d_{i,j}L$ gluing together along an arc. From Lemma~\ref{lemma4.2} together with Proposition~\ref{proposition2.2},
$$
X^L_{i,j}\simeq \widehat{S^3\smallsetminus d_{i,j}L}
$$
with $\nu(d_{i,j}(L))-\nu(d_i(L))-\nu(d_j(L))+\nu(L)=0$ in this case.

\noindent\textbf{Case 2.} \textit{One and only one of $l_i$ and $l_j$ is splittable in $L$}. We may assume that $l_i$ is splittable in $L$ and $l_j$ is not splittable in $L$. Then there is a connected sum decomposition
\begin{equation}\label{equation4.4}
S^3\smallsetminus L\cong (S^3\smallsetminus d_iL)\# (S^3\smallsetminus l_i)
\end{equation}
with $\nu(L)=\nu(d_iL)+1$. Consider the connected sum decomposition
\begin{equation}\label{equation4.5}
S^3\smallsetminus d_iL \cong (S^3\smallsetminus L^{[1]})\# (S^3\smallsetminus L^{[2]})\#\cdots \# (S^3\smallsetminus L^{[\nu(d_iL)+1]}).
\end{equation}
Then $l_j$ lies in one of the sublinks $L^{[1]},\ldots, L^{[\nu(d_iL)+1]}$. We may assume that $l_j\in L^{[1]}$. This gives a connected sum decomposition
\begin{equation}\label{equation4.6}
S^3\smallsetminus d_{l_j}L^{[1]}\cong (S^3\smallsetminus L^{[1,1]})\#(S^3\smallsetminus L^{[1,2]})\#\cdots\# (S^3\smallsetminus L^{[1,h]})
\end{equation}
with $h\geq 1$. By inputting this decomposition formula into equations ~(\ref{equation4.5}) and ~(\ref{equation4.4}), we have
$$
\nu(d_jL)=h+\nu(d_iL)\geq \nu(L)
$$
and $\nu(d_{i,j}L)=h+\nu(d_iL)-1$. Let $A$ be the homotopy pushout of
$$
\bigvee_{k=1}^{\nu(d_iL)}S^2\lTo^{p} \bigvee_{k=1}^{\nu(L)=\nu(d_iL)+1}S^2\rTo^{f} \bigvee_{k=1}^{\nu(d_{j}L)=h+\nu(d_iL)} S^2
$$
by taking out the $2$-spheres as in the $2$-corner ~(\ref{equation4.3}), where $p$ is a canonical retraction by pinching one of the $2$-sphere to a point and $f$ is a canonical inclusion. Then
$$
A\simeq \bigvee_{k=1}^{\nu(d_{j}L)-1=\nu(d_{i,j}L)}S^2,
$$
which is homotopy equivalent to the union of the separating $2$-spheres in $S^3\smallsetminus d_{i,j}L$ gluing together along an arc. Thus
$$
X^L_{i,j}\simeq \widehat{S^3\smallsetminus d_{i,j}L}
$$
with
$$
\nu(d_{i,j}(L))-\nu(d_i(L))-\nu(d_j(L))+\nu(L)=0
$$
in this case.

\noindent\textbf{Case 3.} \textit{Both $l_i$ and $l_j$ are not splittable in $L$}. Let
$$
S^3\smallsetminus L\cong (S^3\smallsetminus L^{[1]})\#(S^3\smallsetminus L^{[2]})\#\cdots\#(S^3\smallsetminus L^{[\nu(L)+1]})
$$
be the complete connected sum decomposition. Since both $l_i$ and $l_j$ are not splittable in $L$, we have $\nu(d_iL)\geq \nu(L)$ and $\nu(d_jL)\geq \nu(L)$.
Since $l_i$ lies in one of $L^{[1]},\ldots,L^{[\nu(L)+1]}$, we may assume that $l_i\in L^{[\nu(L)+1]}$ and let
$$
S^3\smallsetminus d_{l_i}L^{[1]}\cong (S^3\smallsetminus \bar L^{[\nu(L)+1]})\#(S^3\smallsetminus \bar L^{[\nu(L)+2]})\#\cdots\#(S^3\smallsetminus \bar L^{[\nu(L)+h+1]})
$$
be the complete connected sum decomposition with $h\geq 0$. Then $$\nu(d_iL)=\nu(L)+h$$ with a connected sum decomposition
\begin{equation}\label{equation4.7}
S^3\smallsetminus d_iL\cong (S^3\smallsetminus \bar L^{[1]})\#(S^3\smallsetminus \bar L^{[2]})\#\cdots\#(S^3\smallsetminus \bar L^{[\nu(L)+h+1]}),
\end{equation}
where $\bar L^{[k]}=L^{[k]}$ for $1\leq k\leq \nu(L)$. Note that
$$
l_j\in \bigcup_{k=1}^{\nu(L)+h+1}\bar L^{[k]}.
$$

\noindent\textbf{Sub-Case 3.1.} $l_j\in \bigcup\limits_{k=1}^{\nu(L)}\bar L^{[k]}$. \textit{In this case $l_i$ and $l_j$ lie in different nonsplittable components of $L$}. We may assume that $l_j\in \bar L^{[1]}=L^{[1]}$. Let
\begin{equation}\label{equation4.8}
S^3\smallsetminus d_{l_j}\bar L^{[1]}\cong (S^3\smallsetminus \tilde L^{-t+1})\# (S^3\smallsetminus \tilde L^{-t+2})\#\cdots\# (S^3\smallsetminus \tilde L^{0})\# (S^3\smallsetminus \tilde L^{1})
\end{equation}
be the complete connected sum decomposition. By inputting decomposition~(\ref{equation4.8}) into decomposition~(\ref{equation4.7}), we have the complete connected sum decomposition
\begin{equation}\label{equation4.9}
S^3\smallsetminus d_{i,j}L\cong (S^3\smallsetminus \tilde L^{-t+1})\# (S^3\smallsetminus \tilde L^{-t+2})\#\cdots\# (S^3\smallsetminus \tilde L^{\nu(L)+h+1})
\end{equation}
with $t\geq 0$, where $\tilde L^{[k]}=\bar L^{[k]}$ for $2\leq k\leq \nu(L)+h+1$. Thus $\nu(d_{i,j}L)=\nu(L)+h+t$. On the other hand, by filling back $l_i$ to $d_{i,j}L$, we have the complete connected sum decomposition
\begin{equation}\label{equation4.10}
S^3\smallsetminus d_jL\cong (S^3\smallsetminus \tilde L^{-t+1})\# (S^3\smallsetminus \tilde L^{-t+1})\#\cdots\#( S^3\smallsetminus \tilde L^{[\nu(L)]})\# (S^3\smallsetminus L^{[\nu(L)+1]}.
\end{equation}
Thus $\nu(d_jL)=\nu(L)+t$.  Let $A$ be the homotopy pushout of
$$
\bigvee_{k=1}^{\nu(d_iL)=\nu(L)+h}S^2\lTo^{f_1} \bigvee_{k=1}^{\nu(L)}S^2\rTo^{f_2} \bigvee_{k=1}^{\nu(d_{j}L)=\nu(L)+t} S^2
$$
by taking out the $2$-spheres as in the $2$-corner ~(\ref{equation4.3}), where both $f_1$ and $f_2$ are the canonical inclusions. Then
$$
A\simeq \bigvee_{k=1}^{\nu(L)+h+t=\nu(d_{i,j}L)}S^2
$$
which is homotopy equivalent to the union of the separating $2$-spheres in $S^3\smallsetminus d_{i,j}L$ gluing together along an arc. Thus
$$
X^L_{i,j}\simeq \widehat{S^3\smallsetminus d_{i,j}L}
$$
with
$$
\nu(d_{i,j}(L))-\nu(d_i(L))-\nu(d_j(L))+\nu(L)=0
$$
in this case.

\noindent\textbf{Sub-Case 3.2.} $l_j\in \bigcup\limits_{k=\nu(L)+1}^{\nu(L)+h+1}\bar L^{[k]}$. \textit{In this case $l_i$ and $l_j$ lie in the same nonsplittable components of $L$}. We may assume that $l_j\in L^{[\nu(L)+h+1]}$. Let
$$
S^3\smallsetminus d_{l_j}\bar L^{[\nu(L)+h+1]}\cong (S^3\smallsetminus \tilde L^{[\nu(L)+h+1]})\# (S^3\smallsetminus \tilde L^{[\nu(L)+h+2]})\#\cdots\# (S^3\smallsetminus \tilde L^{[\nu(L)+h+t+1]})
$$
be the complete connected sum decomposition. Similar to the arguments in the above case, we have
$$
\nu(d_{i,j}L)=\nu(L)+h+t.
$$
Consider the complete connected sum decomposition
$$
S^3\smallsetminus d_{l_i,l_j} L^{[\nu(L)+1]}\cong (S^3\smallsetminus \tilde L^{[\nu(L)+1]})\# (S^3\smallsetminus \tilde L^{[\nu(L)+2]})\#\cdots\#(S^3\smallsetminus \tilde L^{[\nu(L)+h+t+1]}),
$$
where $\tilde L^{[k]}=\bar L^{[k]}$ for $\nu(L)+1\leq k\leq \nu(L)+h$. Observe that $l_j$ is linked with the nonsplittable components $\tilde L^{[\nu(L)+h+1]},\ldots,\tilde L^{[\nu(L)+h+t+1]}$. For having nonsplittable property of $L^{[\nu(L)+1]}$, $l_i$ must be linked with each of the remaining nonsplittable components
$$
\tilde L^{[\nu(L)+1]}=\bar L^{[\nu(L)+1]},\ldots,\tilde L^{[\nu(L)+h]}=\bar L^{[\nu(L)+h]}
$$
as well as the union
\begin{equation}\label{equation4.11}
l_j\cup \bigcup_{k=\nu(L)+h+1}^{\nu(L)+h+t+1}\tilde L^{[k]}.
\end{equation}

\noindent\textbf{Sub-Sub-Case 3.2.1.} $l_i$ can be deformed away from the sublink $$\bigcup_{k=\nu(L)+h+1}^{\nu(L)+h+t+1}\tilde L^{[k]}.$$ In this case, \textit{ $l_i$ and $l_j$ are linked together lying in the same nonsplittable components of $L$ with the property that none of nonsplittable components of $d_{i,j}L$ can be linked by both $l_i$ and $l_j$.} Then, by filling back $l_i$, we have the completed connected sum decomposition
$$
S^3\smallsetminus d_{l_j} L^{[\nu(L)+1]}=(S^3\smallsetminus L')\# (S^3\smallsetminus \tilde L^{[\nu(L)+h+1]})\#\cdots \#(S^3\smallsetminus \tilde L^{[\nu(L)+h+t+1]}),
$$
where
$$
L'=l_i\cup \bigcup_{k=\nu(L)+1}^{\nu(L)+h}\tilde L^{[k]}
$$
which is not splittable. Thus $\nu(d_j(L))=\nu(L)+t+1$.  Let $A$ be the homotopy pushout of
$$
\bigvee_{k=1}^{\nu(d_iL)=\nu(L)+h}S^2\lTo^{g_1} \bigvee_{k=1}^{\nu(L)}S^2\rTo^{g_2} \bigvee_{k=1}^{\nu(d_{j}L)=\nu(L)+t+1} S^2
$$
by taking out the $2$-spheres as in the $2$-corner ~(\ref{equation4.3}), where $g_i$ are canonical inclusions. Then
$$
A\simeq \bigvee_{k=1}^{\nu(L)+h+t+1}.
$$
From Lemma~\ref{lemma4.2}, we have
$$
X^L_{i,j}\simeq \widehat{S^3\smallsetminus d_{i,j}L}\vee S^3
$$
with $$
\nu(d_{i,j}(L))-\nu(d_i(L))-\nu(d_j(L))+\nu(L)=-1
$$
in this case.

\noindent\textbf{Sub-Sub-Case 3.2.2.} \textit{$l_i$ can NOT be deformed away from the sublink $$\bigcup_{k=\nu(L)+h+1}^{\nu(L)+h+t+1}\tilde L^{[k]}$$
and is linked with only one of the nonsplittable components $$\tilde L^{[\nu(L)+h+1]},\ldots,\tilde L^{[\nu(L)+h+t+1]}.$$} Then, similar to the above arguments, we have  $\nu(d_j(L))=\nu(L)+t$ and
$$
X^L_{i,j}\simeq \widehat{S^3\smallsetminus d_{i,j}L}
$$
with $\nu(d_{i,j}L)-\nu(d_iL)-\nu(d_jL)+\nu(L)=0$ in this case.

\noindent\textbf{Sub-Sub-Case 3.2.3.} $l_i$ can NOT be deformed away from the sublink $$\bigcup_{k=\nu(L)+h+1}^{\nu(L)+h+t+1}\tilde L^{[k]}$$
and is linked with more than one of the nonsplittable components $$\tilde L^{[\nu(L)+h+1]},\ldots,\tilde L^{[\nu(L)+h+t+1]}.$$
 In this case, \textit{both $l_i$ and $l_j$ lie in the same nonsplittable component of $L$ and there are at least two nonsplittable components in $d_{i,j}L$ that are linked by both $l_i$ and $l_j$.} Then, similar to the above arguments, we have  $\nu(d_j(L))<\nu(L)+t$ and
$$
X^L_{i,j}\simeq \widehat{S^3\smallsetminus d_{i,j}L}\vee \bigvee_{k=1}^{\nu(d_{i,j}(L))-\nu(d_i(L))-\nu(d_j(L))+\nu(L)}S^2
$$
with $\nu(d_{i,j}L)-\nu(d_iL)-\nu(d_jL)+\nu(L)>0$ in this case.
\end{proof}

\begin{thm}\label{theorem4.4}
Let $L=\{l_1,\ldots,l_n\}$ be an $n$-link in $S^3$ with $n\geq 2$. Let $$L'=L\smallsetminus \{l_i,l_j\}$$ with $1\leq i\not=j\leq n$.
\begin{enumerate}
\item If $n\leq 3$, then $\calA(L,L')=0$ and so $A(L,l_i)\cap A(L,l_j)=[A(L,l_i),A(L,l_j)]$.
\item Let $n>3$. Then the following statements are equivalent:
\begin{enumerate}
\item[(i)] $\calA(L,L')\not=0$.
\item[(ii)] $\nu(L')-\nu(L\smallsetminus\{l_i\})-\nu(L\smallsetminus\{l_j\})+\nu(L)>0$.
\item[(iii)] $l_i$ and $l_j$ are linked together lying in the same nonsplittable components of $L$ with the property that there are two or more nonsplittable components of $d_{i,j}L$ that can be linked with both $l_i$ and $l_j$.
\item[(iv)] $\calA(L,L')$ is isomorphic to
$$\pi_2\left(\bigvee_{k=1}^{\nu(d_{i,j}(L))-\nu(d_i(L))-\nu(d_j(L))+\nu(L)}G(L')\ltimes S^2\right),$$
which is a free abelian group of (countably) infinite rank.
\end{enumerate}
\end{enumerate}
\end{thm}
\begin{proof}
By~\cite[Corollary 3.4]{brownloday}, $\calA(L,L')\cong \pi_2(X^L_{i,j})$. The assertion follows from Lemma~\ref{lemma4.3}.
\end{proof}

\begin{cor}
Let $L=\{l_1,\ldots,l_n\}$ be an $n$-link in $S^3$ with $n\geq 2$. Let $$L'=L\smallsetminus \{l_i,l_j\}$$ with $1\leq i\not=j\leq n$. Then
the quotient group
$$A(L,l_i)\cap A(L,l_j)/[A(L,l_i),A(L,l_j)]$$
is either trivial or a free abelian group of (countably) infinite rank.
\end{cor}

\subsection{The Groups $\calA(L,L')$ for Sub $3$-Links $L'$}
Let $L$ be an $n$-link in $S^3$ with $n\geq3$ and let $\{i,j,k\}\subseteq \{1,2,\ldots,n\}$ be three distinct labels. Let $X^L_{i,j,k}$ be the homotopy pushout of the $3$-corner
\begin{diagram}
\widehat{S^3\smallsetminus L}&     & \rTo   &     & \widehat{S^3\smallsetminus d_iL}&    &\\
        &\rdTo&         &     &\dDashto      &\rdTo&\\
\dTo    &     &\widehat{S^3\smallsetminus d_jL}&\rTo &          &       &\widehat{S^3\smallsetminus d_{i,j}L}\\
        &     &     &     &     &       &\\
\widehat{S^3\smallsetminus d_kL}& \rDashto   &     &     &\widehat{S^3\smallsetminus d_{i,k}L}&&\\
         &\rdTo& \dTo    & &&&\\
         &     &\widehat{S^3\smallsetminus d_{j,k}L}.& &&&\\
\end{diagram}
Our first step is to determine the homotopy type of $X^L_{i,j,k}$. From Lemma~\ref{lemma4.1}, there is a homotopy pushout diagram
\begin{equation}\label{equation4.12}
\begin{diagram}
X^{L}_{i,j}&\rTo&\widehat{S^3\smallsetminus d_{i,j}L}\\
\dTo&&\dTo\\
X^{d_kL}_{i,j}&\rTo&X^L_{i,j,k}.\\
\end{diagram}
\end{equation}
Lemma~\ref{lemma4.3} can be reformulated as follows: Define $\nu(\emptyset)=-1$ for an empty link. Let
$$
\chi^L_{i,j}=\nu(d_{i,j}L)-\nu(d_iL)-\nu(d_jL)+\nu(L).
$$
Then, by Lemma~\ref{lemma4.3}, we have
\begin{equation}\label{equation4.13}
X^L_{i,j}\simeq \widehat{S^3\smallsetminus d_{i,j}L}\vee\bigvee_{s=1}^{|\chi^L_{i,j}|} S^{2+\frac{|\chi^L_{i,j}|-\chi^L_{i,j}}{2}}
\end{equation}
for any $n$-link $L$ with $n\geq 2$ and any $\{i,j\}\subseteq\{1,2,\ldots,n\}$ with $i\not=j$. Let
\begin{equation}\label{equation4.14}
\begin{array}{rcl}
\chi^L_{i,j,k}&=&\chi^{d_kL}_{i,j}-\chi^L_{i,j}\\
&=&\nu(d_{i,j,k}L)-\nu(d_{i,j}L)-\nu(d_{i,k}L)-\nu(d_{j,k}L)\\
& &\ +\nu(d_iL)+\nu(d_jL)+\nu(d_kL)-\nu(L).\\
\end{array}
\end{equation}

\begin{prop}\label{proposition4.6}
Let $L$ be an $n$-link in $S^3$ with $n\geq3$ and let $\{i,j,k\}\subseteq \{1,2,\ldots,n\}$ be three distinct labels.
\begin{enumerate}
\item If $\chi^L_{i,j}=\nu(L)-\nu(d_iL)-\nu(d_jL)+\nu(L)=0$, then
$$
X^L_{i,j,k}\simeq \widehat{S^3\smallsetminus d_{i,j,k}L}\vee \bigvee_{s=1}^{|\chi^{d_L}_{i,j}|} S^{2+\frac{|\chi^{d_kL}_{i,j}|-\chi^{d_kL}_{i,j}}{2}}
$$
with $\chi^L_{i,j,k}=\chi^{d_kL}_{i,j}$.
\item If $\chi^L_{i,j}=\nu(L)-\nu(d_iL)-\nu(d_jL)+\nu(L)=-1$, then
$$
X^L_{i,j,k}\simeq \widehat{S^3\smallsetminus d_{i,j,k}L} \simeq K( G(d_{i,j,k}L),1)
$$
with $\chi^L_{i,i,k}=0$.
\end{enumerate}
\end{prop}
\begin{proof}
(1). By Lemma~\ref{lemma4.3}(3iii), the top map in the push-out diagram~(\ref{equation4.12}) is a homotopy equivalence. Thus
$$
X^L_{i,j,k}\simeq X^{d_kL}_{i,j}
$$
and the assertion follows from equation~(\ref{equation4.13}).

(2). We first show that $\chi^{d_kL}_{i,j}=-1$. According to Lemma~\ref{lemma4.3}(3i), there exists a nonsplittable component $\bar L$ of $L$ with property that $\{l_i,l_j\}\subseteq \bar L$, $l_i$ is linked with $l_j$ and there are no nonsplittable components in $d_{l_i,l_j}\bar L$ that are linked by both $l_i$ and $l_j$. If $l_k\in L\smallsetminus \bar L$, then $\bar L$ is also a nonsplittable component of $d_kL$ with the above property. By Lemma~\ref{lemma4.3}(3i), we have $\chi^{d_kL}_{i,j}=-1$. Suppose that $l_k\in \bar L$. Since $l_i$ and $l_j$ are linked, $\{l_i,l_j\}$ must be in the same nonsplittable component of $d_{l_k}\bar L$. If there exists a nonsplittable component $L'$ of $d_{l_i,l_j}(d_{l_k}L)$ that could be linked by both $l_i$ and $l_j$, then, by filling back $l_k$, there exists a nonsplittable component $L''$ of $d_{l_i,l_j}\bar L$ with $L'\supseteq L'$ that is linked by both $l_i$ and $l_j$ which contradicts to the above property for $l_i$ and $l_j$. Thus there are no nonsplittable components of $d_{l_i,l_j}(d_{l_k}L)$ that can be linked by both $l_i$ and $l_j$. Hence $\chi^{d_kL}_{i,j}=-1$.

From Lemma~\ref{lemma4.3}(3i) together with diagram~(\ref{equation4.12}), there is a homotopy pushout diagram
\begin{diagram}
X^{L}_{i,j}\simeq \widehat{S^3\smallsetminus d_{i,j}L}\vee S^3 &\rTo&\widehat{S^3\smallsetminus d_{i,j}L}\\
\dTo>{f}&&\dTo\\
X^{d_kL}_{i,j}\simeq \widehat{S^3\smallsetminus d_{i,j,k}L}\vee S^3 &\rTo&X^L_{i,j,k}\\
\end{diagram}
in which $f$ maps $S^3$ to $S^3$ of degree $1$. Thus $X^L_{i,j,k}\simeq \widehat{S^3\smallsetminus d_{i,j,k}L}$ and hence the result.
\end{proof}

Now we consider the case for $\chi^L_{i,j}>0$. By Lemma~\ref{lemma4.3}(3ii), there is a complete splitting decomposition
\begin{equation}\label{equation4.15}
d_{i,j}L\cong L^{[-a+1]}\sqcup\cdots \sqcup L^{[0]}\sqcup \cdots \sqcup L^{[b]}\sqcup\cdots \sqcup L^{[b+c]}\sqcup L^{[b+c+1]}\sqcup\cdots \sqcup L^{[b+c+d]},
\end{equation}
where
\begin{enumerate}
\item each factor $L^{[s]}$ is nonsplittable,
\item $a,c,d\geq 0$ and $b\geq 2$,
\item $l_i$ is linked with each of $L^{[-a+1]},L^{[-a+2]},\ldots, L^{[b]}$,
\item $l_j$ is linked with each of $L^{[1]}, L^{[2]},\ldots, L^{[b+c]}$.
\end{enumerate}
Then
\begin{equation}\label{equation4.16}
\begin{array}{rcl}
\nu(L)&=&d\\
\nu(d_iL)&=&c+d\\
\nu(d_jL)&=&a+d\\
\nu(d_{i,j}L)&=&a+b+c+d-1\\
\end{array}
\end{equation}
with $\chi^L_{i,j}=b-1$. Note that $b$ is the number of the nonsplittable components of $d_{i,j}L$ that are linked with both $l_i$ and $l_j$.

\begin{prop}\label{proposition4.7}
Let $L$ be an $n$-link in $S^3$ with $n\geq3$ and let $\{i,j,k\}\subseteq \{1,2,\ldots,n\}$ be three distinct labels. Suppose that $\chi^L_{i,j}>0$.
\begin{enumerate}
\item $\chi^L_{i,j,k}\geq -1$.
\item If $\chi^L_{i,j,k}=-1$, then $$X^L_{i,j,k}\simeq \widehat{S^3\smallsetminus d_{i,j,k}L}\vee S^3.$$
\item If $\chi^L_{i,j,k}\geq 0$, then
$$
X^L_{i,j,k}\simeq \widehat{S^3\smallsetminus d_{i,j,k}L}\vee\bigvee_{s=1}^{\chi^L_{i,j,k}}S^2.
$$
\end{enumerate}
\end{prop}
\begin{proof}
Since $\chi^L_{i,j}=b-1>0$, there are at least two nonsplittable components of $d_{i,j}L$ that are linked with both $l_i$ and $l_j$. By removing $d_k$, there is at least one nonsplittable components of $d_{i,j,k}L$ that are linked with both $l_i$ and $l_j$. Thus $l_i$ and $l_j$ must lie in the same nonsplittable component of $d_kL$.

\noindent\textbf{Case 1.} \textit{$l_k\in L^{[t]}$ for some $t$ with $b+1\leq t\leq b+c+d$ or $-a+1\leq t\leq 0$}. The number of the nonsplittable components of $d_{i,j}(d_kL)=d_{i,j,k}L$ remains the same as $b$. Thus  we have
$$
\chi^{d_kL}_{i,j}=b-1
$$
and so
$$
\chi^L_{i,j,k}=\chi^{d_KL}_{i,j}-\chi^{L}_{i,j}=0
$$
in this case.

\noindent\textbf{Case 2.} \textit{$l_k\in L^{[t]}$ for some $t$ with $1\leq t\leq b$}. We may assume that $l_k\in L^{[1]}$. Let
$$
d_kL^{[1]}=L^{[1,1]}\sqcup L^{[1,2]}\sqcup\cdots\sqcup L^{[1,e+1]}
$$
be the complete splitting decomposition with $f$, $0\leq f\leq e+1$, factors that are linked by both $l_i$ and $l_j$. (\textbf{Note.} If $d_kL^{[1]}=\emptyset$, then $f=0$.) Then the number of nonsplittable components of $d_{i,j,k}L$ linked with both $l_i$ and $l_j$ is $b+f-1$. Thus
$$
\chi^{d_kL}_{i,j}=b+f-1
$$
and so
$$
\chi^L_{i,j,k}=\chi^{d_kL}_{i,j}-\chi^{L}_{i,j}=f-1.
$$

From the both cases, we have $\chi^L_{i,j,k}\geq-1$. This proves assertion (1). Assertions (2) and (3) then follow from diagram~(\ref{equation4.12}) and Lemma~\ref{lemma4.3}.
\end{proof}

By putting Propositions~\ref{proposition4.6} and~\ref{proposition4.7} together, we have the following.

\begin{prop}\label{proposition4.8}
Let $L$ be an $n$-link in $S^3$ with $n\geq3$. Then
$$
X^L_{i,j,k}\simeq \widehat{S^3\smallsetminus d_{i,j,k}L}\vee \bigvee_{s=1}^{|\chi^L_{i,j,k}|} S^{2+\frac{|\chi^L_{i,j,k}|-\chi^L_{i,j,k}}{2}}.
$$
for any distinct labels $\{i,j,k\}\subseteq \{1,2,\ldots,n\}$.\hfill $\Box$
\end{prop}

Recall that $A(L,l_i)$ is defined to be the kernel of $G(L)\to G(d_iL)$. Denote $A(L,l_i)$ by $A_i$. Let
$$
\bar\calA(L, L\smallsetminus d_{i,j,k}L)=(A_i\cap A_j\cap A_k)/\left([A_i\cap A_j,A_k]\cdot [A_i\cap A_k,A_j]\cdot [A_j\cap A_k,A_i]\right).
$$
Observe that
$$[[A_i,A_j,A_k]_S=[[A_i,A_j],A_k]\cdot [[A_i,A_k],A_j]\cdot [[A_j,A_k],A_i]$$ is a (normal) subgroup of
$$
[A_i\cap A_j,A_k]\cdot [A_i\cap A_k,A_j]\cdot [A_j\cap A_k,A_i].
$$
Thus there is a canonical epimorphism
$$
\calA(L,L\smallsetminus d_{i,j,k}L)\twoheadrightarrow\bar \calA(L,L\smallsetminus d_{i,j,k}L)
$$
with the kernel given by
$$
\left([A_i\cap A_j,A_k]\cdot [A_i\cap A_k,A_j]\cdot [A_j\cap A_k,A_i]\right)/[[A_i,A_j,A_k]_S.
$$

\begin{thm}\label{theorem4.9}
Let $L$ be an $n$-link in $S^3$ with $n\geq3$ and let $\{i,j,k\}\subseteq \{1,2,\ldots,n\}$ be three distinct labels. Then there is an isomorphism of $\Z[G(L)]$-modules
$$
 \bar\calA(L, L\smallsetminus d_{i,j,k}L)\cong \pi_3\left(\bigvee_{s=1}^{|\chi^L_{i,j,k}|} G(d_{i,j,k}L)\ltimes S^{2+\frac{|\chi^L_{i,j,k}|-\chi^L_{i,j,k}}{2}}\right).
$$
In particular, $\bar\calA(L, L\smallsetminus d_{i,j,k}L)=0$ if and only if $\chi^L_{i,j,k}=0$.
\end{thm}
\begin{proof}
The assertion follows from Proposition~\ref{proposition4.8} and~\cite[Theorem 1]{EM}.
\end{proof}

\begin{cor}
Let $L$ be an $n$-link in $S^3$ with $n\geq3$ and let $\{i,j,k\}\subseteq \{1,2,\ldots,n\}$ be three distinct labels. Then $\bar\calA(L, L\smallsetminus d_{i,j,k}L)$ is either $0$, $\Z$ or a free abelian group of (countably) infinite rank.\hfill $\Box$
\end{cor}

\begin{example}
{\rm In this example, we determine the group $\calA(L,L)$ for any $3$-link $L$.

Let $L$ be a $3$-link and let $(i,j,k)=(1,2,3)$. Then $d_{1,2,3}L=\emptyset$ and so
$$
G(d_{1,2,3}L)=\pi_1(S^3)=0.
$$
It follows that
$$
\bar\calA(L, L)\cong \pi_3\left(\bigvee_{s=1}^{|\chi^L_{1,2,3}|} S^{2+\frac{|\chi^L_{1,2,3}|-\chi^L_{1,2,3}}{2}}\right).
$$
By Theorem~\ref{theorem4.4}(1), we have
$$
\calA(L,L\smallsetminus d_{i,j}L)=0
$$
for $1\leq i\not=j\leq 3$. Thus
$$
A_i\cap A_j=[A_i,A_j]
$$
for $1\leq i\not=j\leq 3$ and so
$$
[A_1\cap A_2,A_3]\cdot [A_1\cap A_3,A_2]\cdot [A_2\cap A_3,A_1]=[[A_1,A_2],A_3]_S.
$$
Thus
$$
\calA(L,L)=\bar\calA(L,L)\cong \pi_3\left(\bigvee_{s=1}^{|\chi^L_{1,2,3}|} S^{2+\frac{|\chi^L_{1,2,3}|-\chi^L_{1,2,3}}{2}}\right).
$$
Note that
$$
\begin{array}{rcl}
\chi^L_{1,2,3}&=&\nu(d_{1,2,3}L)-\nu(d_{1,2}L)-\nu(d_{1,3}L)-\nu(d_{2,3}L)\\
&& \ +\nu(d_1L)+\nu(d_2L)+\nu(d_3L)-\nu(L)\\
&=&-1-0-0-0+\nu(d_1L)+\nu(d_2L)+\nu(d_3L)-\nu(L)\\
\end{array}
$$
because $d_{1,2,3}L=\emptyset$ and $d_{i,j}L$ is a knot.

\noindent\textbf{Case 1.} $\nu(L)=2$. In this case $L$ is given by three knots with none of them linked each other. Then
$$
\nu(d_1L)=\nu(d_2L)=\nu(d_3L)=1
$$
and so $\chi^L_{1,2,3}=0$. By Theorem~\ref{theorem4.9}, $\calA(L,L)=\bar\calA(L,L)=0$ and so
$$
A_1\cap A_2\cap A_3=[[A_1,A_2],A_3]_S.
$$

\noindent\textbf{Case 2.} $\nu(L)=1$. In this case $L$ consists of one knot together with two other linked knots. We may assume that $l_1$ and $l_2$ are linked and $l_3$ can be deformed away from $\{l_1,l_2\}$. Then $\nu(d_1L)=\nu(d_2L)=1$ and $\nu(d_3L)=0$ and so
$$
\chi^L_{1,2,3}=-1+1+1-1=0.
$$
Similar to Case 1, we have
$$
A_1\cap A_2\cap A_3=[[A_1,A_2],A_3]_S.
$$

\noindent\textbf{Case 3.} $\nu(L)=0$. In this case $L$ is not a splittable $3$-link. Then
$$
-1\leq \chi^L_{1,2,3}=\nu(d_1L)+\nu(d_2L)+\nu(d_3L)-1\leq 2.
$$

If $\chi^L_{1,2,3}=1$, then it requires that two of $\nu(d_1L), \nu(d_2L), \nu(d_3L)$ are $1$ and one of them is $0$. We may assume that $\nu(d_1L)=\nu(d_2L)=1$ and $\nu(d_3L)=0$. From $\nu(d_1L)=1$, $l_2$ and $l_3$ are not linked each other. Since $L$ is not splittable, then $l_1$ linked with both $l_2$ and $l_3$. By removing $l_2$, we have $\nu(d_2L)=0$ because $l_1$ and $l_3$ are linked. Thus $\chi^L_{1,2,3}\not=1$.

\noindent\textbf{Sub-Case 3.1.} $\chi^L_{1,2,3}=2$. This equality holds if and only if $$\nu(d_1L)=\nu(d_2L)=\nu(d_3L)=1.$$ In other words, $L$ is a nonsplittable $3$-link such that it becomes splittable by removing any one of its components. The Brunnian $3$-links are the examples in this case. By Theorem~\ref{theorem4.9}, we have
$$
\begin{array}{rcl}
\calA(L,L)&=&\bar\calA(L,L)\\
&\cong& \pi_3(S^2\vee S^2)\\
&\cong& \pi_2(\Omega(S^2\vee S^2))\\
&\cong &\pi_2(\Omega S^2)\oplus \pi_2(\Omega S^2)\oplus \pi_2(\Omega\Sigma (\Omega S^2) \wedge(\Omega S^2))\\
&=&\Z\oplus\Z\oplus\Z.\\
\end{array}
$$

\noindent\textbf{Sub-Case 3.2.} $\chi^L_{1,2,3}=0$. In this case, it requires that one of $\nu(d_1L), \nu(d_2L), \nu(d_3L)$ are $1$ and two of them is $0$. In other words, $L$ is a nonsplittable $3$-link such that it becomes splittable by removing one component but it stands nonsplittable by removing any of the other two components. In this case, we have $$A_1\cap A_2\cap A_3=[[A_1,A_2],A_3]_S.$$

\noindent\textbf{Sub-Case 3.3.} $\chi^L_{1,2,3}=-1$. This equality holds if and only if $$\nu(d_1L)=\nu(d_2L)=\nu(d_3L)=0,$$ if and only if $L$ is a strongly nonsplittable $3$-link. The Hopf $3$-link, which is given by taking the pre-image of the $3$ points in $S^2$ from the Hopf map $S^3\to S^2$, is an example in this case. In this case,
$$
\calA(L,L)\cong \pi_3(S^3)=\Z
$$
either by Theorem~\ref{theorem1.1} or ~\ref{theorem4.9}. The computation is finished.\hfill $\Box$
}\end{example}

\section{Homotopy-group Invariants and Milnor's Invariants}\label{section5}
Let $L$ be an $n$-link in $S^3$. Recall that Milnor's link group~\cite{Milnor}, denoted by $\calG(L)$, is defined as follows: Let $A_i=A(L,l_i)$ be the kernel of $G(L)\to G(d_iL)$. Then
$\calG(L)$ is defined by
$$
\calG(L)=G(L)/\prod_{i=1}^n [A_i,A_i].
$$
In this section, we construct examples of $(n+1)$-links given in the form $L\cup l$ that have nontrivial homotopy-group invariants with the property that $l$ represents the trivial element in Milnor's link group $\calG(L)$. In other words, we provide examples of the links $L\cup l$ labeled by the nontrivial elements in $\pi_n(S^3)$ in which $l$ is linked with  $L$ but as homotopy links~\cite{HL} in Milnor's sense $l$ is unlinked with $L$.

Consider the Hopf fibration
\begin{equation}\label{equation5.1}
p\colon S^3\longrightarrow S^2.
\end{equation}
Let
$$
Q_n=\{q_1,\ldots,q_n\}\subseteq S^2
$$
be the $n$ distinct points in $S^2$. Let
\begin{equation}\label{equation5.2}
L_n=p^{-1}(Q_n)
\end{equation}
Then $L_n=\{l_1,\ldots,l_n\}$ is an $n$-link in $S^3$, where $l_i=p^{-1}(q_i)$. Let $I=\{i_1,\ldots,i_k\}$ be any nonempty subset of $\{1,\ldots,n\}$. From the fibration
\begin{equation}\label{equation5.3}
S^1\longrightarrow S^3\smallsetminus p^{-1}(q_{i_1},\ldots,q_{i_k})\longrightarrow S^2\smallsetminus \{q_{i_1},\ldots,q_{i_k}\},
\end{equation}
the space $S^3\smallsetminus p^{-1}(q_{i_1},\ldots,q_{i_k})$ is a $K(\pi,1)$-space. By Theorem~\ref{theorem2.1}, the sublink
$$
\{l_{i_1},\ldots,l_{i_k}\}
$$
is nonsplittable. It follows that $L_n$ is a strongly nonsplittable $n$-link.

\begin{thm}\label{theorem5.1}
Let $\phi\colon G(L_n)\to \calG(L_n)$ be the quotient homomorphism Then
$$
A_1\cap A_2\cap\cdots \cap A_n\subseteq \Ker(\phi)
$$
for $n\geq 4$. Thus Milnor's link group $\calG(L_n)$ gives no information for
$$
\calA(L_n,L_n)\cong \pi_n(S^3)
$$
when $n\geq4$.
\end{thm}
For a group $G$, let $\gamma^1(G)=G$ and $\gamma^{i+1}(G)=[\gamma^i(G),G]$.
\begin{proof}
From fibration~(\ref{equation5.3}), there is a commutative diagram of short exact sequences of groups
\begin{diagram}
\pi_1(S^1)&\rInto& \pi_1(S^3\smallsetminus L_n)&\rOnto& \pi_1(S^2\smallsetminus Q_n)\\
\dEq&&\dTo & &\dTo\\
\pi_1(S^1)&\rInto& \pi_1(S^3\smallsetminus p^{-1}(q_1))&\rOnto& \pi_1(S^2\smallsetminus\{q_1\})=\{1\}.\\
\end{diagram}
It follows that the link group
\begin{equation}\label{equation5.4}
G(L_n)=\pi_1(S^3\smallsetminus L_n)\cong \pi_1(S^1)\times \pi_1(S^2\smallsetminus Q_n)\cong \Z\times \pi_1(S^2\smallsetminus Q_n).
\end{equation}
Note that
$$
\pi_1(S^2\smallsetminus Q_n)\cong \pi_1(\R^2\smallsetminus Q_{n-1})
$$
is a free group of rank $n-1$. By~\cite[Theorem 5.7]{MKS},
$$
\gamma^{n-1}\pi_1(S^2\smallsetminus Q_n)/\gamma^n\pi_1(S^2\smallsetminus Q_n)
$$
is a free abelian group. From decomposition~(\ref{equation5.4}),
\begin{equation}\label{equation5.5}
\gamma^{n-1}G(L_n)/\gamma^nG(L_n)\cong \gamma^{n-1}\pi_1(S^2\smallsetminus Q_n)/\gamma^n\pi_1(S^2\smallsetminus Q_n) \textrm{ is torsion free.}
\end{equation}

We claim that
\begin{equation}\label{equation5.6}
A_1\cap\cdots\cap A_n\leq \gamma^nG(L_n).
\end{equation}
From Theorem~\ref{theorem1.1}(1),
$$
A_1\cap\cdots\cap A_n\leq A_1\cap\cdots\cap A_{n-1}=[[A_1,A_2],\ldots,A_{n-1}]_S\leq \gamma^{n-1}G(L_n).
$$
Note that $A_1\cap\cdots\cap A_n\leq \gamma^nG(L_n)$ if and only if $A_1\cap\cdots\cap A_n$ has the trivial image in the quotient group
$
\gamma^{n-1}G(L_n)/\gamma^nG(L_n).
$
Suppose that there exists an element $\alpha\in A_1\cap\cdots\cap A_n$ that has nontrivial image in $\gamma^{n-1}G(L_n)/\gamma^nG(L_n)$. Then, for each $k>0$,
$
\alpha^k
$
has nontrivial image in $\gamma^{n-1}G(L_n)/\gamma^nG(L_n)$ because $\gamma^{n-1}G(L_n)/\gamma^nG(L_n)$ is torsion free by~(\ref{equation5.5}). Note that
$$
[[A_1,A_2],\ldots,A_n]_S\leq \gamma^nG(L_n).
$$
Thus
$$
\alpha^k\not\in [[A_1,A_2],\ldots,A_n]_S
$$
and so the quotient group
$$
\calA(L_n,L_n)=A_1\cap\cdots\cap A_n/[[A_1,A_2],\ldots,A_n]_S
$$
contains an element of infinite order. By Theorem~\ref{theorem1.1}(2), $\calA(L_n,L_n)\cong \pi_n(S^3)$. By the Serre Theorem, $\pi_n(S^3)$ is a finite group. This gives a contradiction. Hence the statement in (\ref{equation5.6}) holds.

Next we claim that
$$
\gamma^n\calG(L_n)=\{1\}.
$$
If this is true, then the assertion follows because $\phi$ maps $\gamma^nG(L_n)$ into $\gamma^n\calG(L_n)$.

Choose the points $q_1,\ldots,q_n$ in the real line $$\R^1\subseteq \R^2\subseteq S^2=\R^2\cup\{\infty\}$$
with $q_1<q_2<\cdots<q_n$. The fundamental group $\pi_1(S^2\smallsetminus Q_n)$ admits a presentation
\begin{equation}\label{equation5.7}
\pi_1(S^2\smallsetminus Q_n)=\la x_1,\ldots,x_n \ | \ x_1x_2\cdots x_n=1\ra,
\end{equation}
where $x_i$ is represented by a loop $\lambda_i$ that goes through the geodesic curve from $\infty$ to a point $q'_i$ in a small neighborhood of $q_i$ followed by a small circle around the point $q_i$ counterclockwise and then return back to $\infty$ through the geodesic curve. Consider the commutative diagram of short exact sequences
\begin{equation}\label{equation5.8}
\begin{diagram}
& &\pi_1(S^1)&\rEq&\pi_1(S^1)\\
&&\dInto&&\dInto\\
A_i&\rInto&G(L_n)&\rOnto& G(d_iL_n)\\
\dTo&&\dOnto>{p_*}&&\dOnto>{p_*}\\
\la\la x_i\ra\ra&\rInto& \pi_1(S^2\smallsetminus Q_n)&\rOnto&\pi_1(S^2\smallsetminus\{q_1,\ldots,q_{i-1},q_{i+1},\ldots,q_n\}),\\
\end{diagram}
\end{equation}
where $\la\la x_i\ra\ra$ is the normal closure of $x_i$ in $\pi_1(S^2\smallsetminus Q_n)$. The epimorphism
$$
p_*\colon G(L_n)\longrightarrow \pi_1(S^2\smallsetminus Q_n)
$$
induces an isomorphism
\begin{equation}\label{equation5.9}
p_*|_{A_i}\colon A_i\rTo^{\cong} \la\la x_i\ra\ra.
\end{equation}
Decomposition~(\ref{equation5.4}) induces a decomposition
$$
\calG(L_n)\cong \Z\times \pi_1(S^2\smallsetminus Q_n)/\prod_{i=1}^n[\la\la x_i\ra\ra,\la\la x_i\ra\ra].
$$
Thus
$$
\gamma^n\calG(L_n)\cong \gamma^n\left(\pi_1(S^2\smallsetminus Q_n)/\prod_{i=1}^n[\la\la x_i\ra\ra,\la\la x_i\ra\ra]\right).
$$
Let
$$
\theta\colon F_{n-1}=\la y_1,\ldots,y_{n-1}\ra\longrightarrow \pi_1(S^2\smallsetminus Q_n)
$$
be the group homomorphism such that $\theta(y_i)=x_i$ for $1\leq i\leq n-1$. From the presentation of $\pi_1(S^2\smallsetminus Q_n)$,  $\theta$ is an isomorphism and so it induces an epimorphism
$$
\bar\theta\colon K_{n-1}=F_{n-1}/\prod_{i=1}^{n-1}[\la\la y_i\ra\ra, \la\la y_i\ra\ra] \rOnto \pi_1(S^2\smallsetminus Q_n)/\prod_{i=1}^n[\la\la x_i\ra\ra,\la\la x_i\ra\ra].
$$
Note that the group $K_{n-1}$ is Milnor's link group of the trivial $(n-1)$-link. According to~\cite[Lemma 5]{Milnor},
$$
\gamma^nK_{n-1}=\{1\}.
$$
From the epimorphism $\bar\theta$, we have
$$
\gamma^n\left(\pi_1(S^2\smallsetminus Q_n)/\prod_{i=1}^n[\la\la x_i\ra\ra,\la\la x_i\ra\ra]\right)=\{1\}
$$
and so is $\gamma^nG(L_n)$. The proof is finished now.
\end{proof}

\begin{example}\label{example5.1}
{\rm
We check the missing case $n=3$ in the above theorem. From the above proof, there is an epimorphism
$$
\bar\theta\colon K_2\rOnto \pi_1(S^2\smallsetminus Q_3)/\prod_{i=1}^3[\la\la x_i\ra\ra,\la\la x_i\ra\ra].
$$
From the presentation of $\pi_1(S^2\smallsetminus Q_3)$, the group $\pi_1(S^2\smallsetminus Q_3)/\prod_{i=1}^3[\la\la x_i\ra\ra,\la\la x_i\ra\ra]$ is obtained by adding the relation
$$
[y_1y_2, w(y_1y_2)w^{-1}]\equiv1
$$
to $K_2$ for any word $w\in K_2$. From the Witt-Hall identities~\cite[Theorem 5.1]{MKS}, since $\gamma_3K_2=\{1\}$,
$$
[y_1y_2, w(y_1y_2)w^{-1}]=[y_1y_2,w][y_1y_2,y_1y_2][y_1y_2,w]^{-1}=1
$$
in the group $K_2$. Thus $\bar\theta$ is an isomorphism in this case. Consider the element
$$
\alpha=[x_1,x_2]\in \pi_1(S^2\smallsetminus Q_3).
$$
Then
$$
\alpha\in \bigcap_{i=1}^3 \la\la x_i\ra\ra=\bigcap_{i=1}^3 A_i
$$
which represents a generator for
$$
\calA(L_3,L_3)\cong \pi_3(S^3).
$$
Moreover the image of $\alpha$ is given by $\bar\theta([y_1,y_2])$. Since $[y_1,y_2]$ is a generator for  $\gamma^2K_2\cong\Z$, the elements
$\alpha^k$ has nontrivial image in $\calG(L_3)$ for $k\not=0$. Hence, in the case $n=3$, the homotopy-group invariants $\pi_3(S^3)$ are detected by Milnor's invariants.\hfill $\Box$
}
\end{example}

\section{Examples of Links Labeled By Homotopy Group Elements and Some Remarks}\label{section6}
\subsection{The Group $\calA(L_n,L_n)$ for $n\leq 5$}
In this subsection, we give some examples of link invariants given by homotopy groups elements. Let $L_n$ be the $n$-link in equation~(\ref{equation5.2}). From Theorem~\ref{theorem1.1}, we have
\begin{equation}\label{equation6.1}
\calA(L_n,L_n)\cong \pi_n(S^3).
\end{equation}
According to~\cite{Toda}, the generator for $\pi_3(S^3)\cong \Z$ is labeled by $\iota$, the generator for $\pi_4(S^3)\cong \Z/2$ is label by $\eta$ and the generator for $\pi_5(S^3)\cong\Z/2$ is labeled $\eta^2$ given by the composite
$$
S^5\rTo^{\eta}S^4\rTo^{\eta}S^3.
$$
Let $l$ be a knot in $S^3\smallsetminus L_n$ that represents an element $\alpha$ in $\calA(L_n,L_n)\cong \pi_n(S^3)$. We label the link $L_n\cup l$ by the homotopy group element $\alpha$. Thus we have $4$-links given in form $L_3\cup l$ labeled by $\iota^k$ for $k\in \Z$. According to Example~\ref{example5.1}, the homotopy group invariants $\iota^k$ give no more information than Milnor's invariants. The first interesting cases are $5$-links labeled by $\eta\in \pi_4(S^3)$ and $6$-links labeled by $\eta^2\in \pi_5(S^3)$. Below we construct the words in $A(L_n,L_n)$, $n=4,5$, in terms of meridians that give the generators for $\calA(L_4,L_4)\cong \calA(L_5,L_5)\cong \Z/2$. Our construction uses the simplicial group techniques. The standard references for the theory of simplicial objects are~\cite{Curtis,May}.

Recall that a \textit{simplicial group} means a sequence of group $\calG=\{G_n\}_{n\geq 0}$ together with \textit{face homomorphisms} $d_i\colon G_n\to G_{n-1}$, $0\leq i\leq n$, and \textit{degeneracy homomorphisms} $s_i\colon G_n\to G_{n+1}$, $0\leq i\leq n$, such that the following \textit{simplicial identities} for $d_i$ and $s_j$ hold. For a simplicial group $\calG$, define the \textit{Moore chains} by
$$
N_n\calG=\bigcap_{i=1}^n\Ker(d_i\colon G_n\to G_{n-1}).
$$
From the simplicial identities, one gets $d_0(N_{n+1}\calG)\unlhd N_n\calG$ for any $n\geq0$ and so a chain complex
\begin{equation}\label{equation6.2}
N\calG=\{N_n\calG\}_{n\geq0}
\end{equation}
of (possibly noncommutative) groups with boundary homomorphism given by $d_0|_{N_n\calG}$. Let
$$
\calZ_n\calG=\bigcap_{i=0}^n\Ker(d_i\colon G_n\to G_{n-1})
$$
and $\calB_n\calG=d_0(N_{n+1}\calG)$. The homology of the chain complex~(\ref{equation6.2}) is given by
$$
H_n(N\calG)=\calZ_n\calG/\calB_n\calG.
$$
The simplicial $n$-simplex with $\Delta[n]$ is defined by
$$
\Delta[n]_k=\{(i_0,i_1,\ldots,i_k)\ | \ 0\leq i_0\leq \cdots\leq i_k\leq n\}
$$
with faces and degeneracies given in the canonical way by removing and doubling the coordinates. Let $\sigma_n=(0,1,\ldots,n)\in \Delta[n]_n$ be the only non-degenerate element in $\Delta[n]$. Let $S^n=\Delta[n]/\partial\Delta[n]$ be the simplicial $n$-sphere. Let $\mathcal{G}=\{G_n\}_{n\geq0}$ be a simplicial group. For any element $z\in G_n$, there is a unique simplicial map
$$
f_z\colon \Delta[n]\longrightarrow \mathcal{G}
$$
such that $f_z(\sigma_n)=z$, called the \textit{representing map for the element $z$}~\cite[Proposition 1.5]{Curtis}. Clearly the representing map $f_z$ factors through the simplicial $n$-sphere $S^n$ if and only if $z\in \mathcal{Z}_n\calG$. Thus each element $z\in \calZ_n\calG$ induces a simplicial map
$$
f_z\colon S^n\longrightarrow \calG
$$
and so a (continuous) map
$$
|f_z|\colon |S^n|\cong S^n\longrightarrow |\calG|
$$
by taking geometric realization. This gives a group homomorphism
\begin{equation}\label{equation6.3}
\calZ_n\calG\longrightarrow \pi_n(|\calG|) \quad z\mapsto [|f_z|].
\end{equation}
The Moore Theorem~\cite{Moore} (also see~\cite[Proposition 5.4]{Kan}) states that the above map induces an isomorphism
\begin{equation}\label{equation6.4}
H_n(N\calG)\cong \pi_n(|\calG|)
\end{equation}
for each $n\geq0$. Namely $\calZ_n\calG\to \pi_n(|\calG|)$ is an epimorphism, and $[|f_z|]=0$ in $\pi_n(|\calG|)$ if and only if $z\in \mathcal{B}_n\calG$.

Note that the geometric realization of a simplicial group is a loop space~\cite[Theorem 3]{Milnor1}. Let $\calG$ be a simplicial group. The $\eta$-operation
$$
\eta^*\colon \pi_n(|\calG|)\longrightarrow \pi_{n+1}(|\calG|)
$$
for $n\geq 1$ is defined as follows: Let $f\colon S^n\to |\calG|\simeq \Omega B|\calG|$ with its adjoint map
$$
f'\colon S^{n+1}\longrightarrow B|\calG|.
$$
Then $\eta^*([f])\in \pi_{n+1}(|\calG|)$ is defined to be the homotopy class represented by the adjoint map of the composite
$$
S^{n+2}\rTo^{\eta}S^{n+1}\rTo^{f'} B|\calG|.
$$

Let $\mathcal{X}$ be a pointed simplicial set. Let $\ast\in X_{0}$ be the
basepoint. The basepoint in $X_{n}$ is $s_{0}^{n}\ast$. Let $F[\mathcal{X}]_{n}$ be the free group generated by $X_{n}$ subject to the single relation that $s_{0}^{n}\ast=1$. (\textbf{Note. }By the simplicial identities, $%
s_{0}^{n}=s_{i_{n}}s_{i_{n-1}}\cdots s_{i_{1}}$ for any sequence $(i_{1},i_{2},\ldots,i_{n})$ with $0\leq i_{k}\leq k-1$.) Then we obtain the simplicial group $F[\mathcal{X}]=\{F[\mathcal{X}]_{n}\}_{n\geq0}$ with the
faces and the degeneracies induced by those of $\mathcal{X}$. The simplicial
group $F[\mathcal{X}]$ is called \textit{Milnor's free group construction}. Clearly the simplicial group $F[\mathcal{X}]$ has the following universal property:
\begin{enumerate}
\item[] \textit{Let $\calG$ be any simplicial group and let $f\colon \mathcal{X}\to \calG$ be any simplicial map such that $f(\ast)=1$. Then there is a unique simplicial homomorphism $\tilde f\colon F[\mathcal{X}]\to G$ with $\tilde f|_{\mathcal{X}}=f$.}
\end{enumerate}
The Milnor~\cite{Milnor2} theorem states that
\begin{equation}\label{equation6.5}
|F[\mathcal{X}]|\simeq \Omega\Sigma |\mathcal{X}|
\end{equation}
for any pointed simplicial set $\mathcal{X}$. (In Milnor's original paper, it requires that $X_0=\{\ast\}$. But this hypothesis can be removed, see for instance~\cite[Theorem 4.9]{Wu}.) The following lemma is useful in our construction.

\begin{lem}\label{lemma6.1}
Let $\calG$ be a simplicial group and let $\alpha$ be an element in $\pi_n(|\calG|)\cong H_n(N\calG)$ represented by an element
$$
z\in \calZ_n\calG.
$$
Suppose that  $n\geq1$. Then the element $[s_0z,s_1z]\in \calZ_{n+1}\calG$ that represents $\eta^*(\alpha)\in \pi_{n+1}(|\calG|)$.
In other words, the operation $\eta^*\colon \pi_n(|\calG|)\to \pi_{n+1}(\calG)$, $\alpha\mapsto \alpha\circ\eta$ is given by the operation $$z\mapsto [s_0z,s_1z]$$ for $z\in \calZ_n\calG$.
\end{lem}

\begin{proof}
Let $f_z\colon S^n\to \mathcal{G}$ be the representing map of the element $z$. Namely $f_z(\bar\sigma_n)=z$, where $\bar\sigma_n$ is the (only) non-degenerate element in $S^n_n$. From the Moore Theorem~(\ref{equation6.4}), the map
$$
|f_z|\colon |S^n|\cong S^n\longrightarrow |\calG|
$$
represents the homotopy class $\alpha$. Let
$$
\tilde f_z\colon F[S^n]\longrightarrow \mathcal{G}
$$
be the (unique) simplicial homomorphism such that $\tilde f_z|_{S^n}=f_z$. According to~\cite[Example 2.21]{Wu1}, the element
$$
[s_0\bar\sigma_n,s_1\bar\sigma_n]\in\calZ_{n+1}F[S^n]
$$
represents the generator $\eta$ for $\pi_{n+1}(F[S^n])$. By applying the simplicial homomorphism $\tilde f_z$, the element
$$
[s_0z,s_1z]=[s_0\tilde f_z(\bar\sigma_n),s_1\tilde f_z(\bar\sigma_n)]=\tilde f_z([s_0\bar \sigma_n,s_1\bar\sigma_n] )
$$
lies in $\mathcal{Z}_{n+1}\calG$ representing the homotopy class $|\tilde f|_*(\eta)\in \pi_{n+1}(|\calG|)$. By~\cite[Proposition 1.1.9]{Wu2}, the map
$$
|f_z|\colon |S^n|\cong S^n\longrightarrow |\mathcal{G}|
$$
extends uniquely (up to homotopy) to an $H$-map
$$
|F[S^n]|\simeq \Omega S^{n+1}\longrightarrow |\mathcal{G}|.
$$
It follows that
$$
|\tilde f_z|\simeq \Omega (|f_z|')\colon \Omega S^{n+1}\longrightarrow |\mathcal{G}|\simeq \Omega B|\mathcal{G}|,
$$
where $|f_z|'\colon S^{n+1}\to B|\mathcal{G}|$ is the adjoint map of $|f_z|$. Thus
$$
|\tilde f_z|_*(\eta)=\Omega (|f_z|')_*(\eta)=\eta^*(\alpha).
$$
The proof is finished.
\end{proof}

Let $\hat F_n$
be the quotient of the free group $F(x_1,\ldots,x_n)$ subject to the single relation
$$x_1x_2\cdots x_n=1.$$ From presentation~(\ref{equation5.7}),  $\pi_1(S^2\smallsetminus Q_n)\cong \hat F_n$. Note that, as a group, $\hat F_n$ is a free group generated by $x_1,\ldots,x_{n-1}$ with additional word $x_n=(x_1\cdots x_{n-1})^{-1}$. From~\cite[subsection 6.1]{BCWW}, the sequence of groups $\hat {\mathcal{F}}=\{\hat F_n\}_{n\geq1}$, with $(\hat{\mathcal{F}})_n=\hat F_{n+1}$, forms a simplicial group in which the faces $d_i\colon \hat F_{n+1}\to \hat F_n$ and the degeneracies $s_i\colon \hat F_{n+1}\to \hat F_{n+2}$ are given
\begin{equation}\label{equation6.6}
\begin{array}{ccc}
d_ix_j=\left\{
\begin{array}{lcl}
x_j&\textrm{ if } &j<i+1,\\
1&\textrm{ if }&j=i+1,\\
x_{j-1}&\textrm{ if }& j>i+1,
\end{array}\right.
& &
s_ix_j=\left\{
\begin{array}{lcl}
x_j&\textrm{ if }& j<i+1,\\
x_jx_{j+1}&\textrm{ if }&j=i+1,\\
x_{j+1}&\textrm{ if }& j>i+1\\
\end{array}\right.\\
\end{array}
\end{equation}
for $0\leq i\leq n$. (\textbf{Note.} In~\cite[subsection 6.1]{BCWW}, the generators for $\hat F_{n+1}$ are labeled by $\hat z_0,\ldots,\hat z_n$. The above formula translates ~\cite[Equation (12)]{BCWW} in terms of our generators $x_1,\ldots,x_{n+1}$.)
 According to~\cite[Proposition 6.1.2]{BCWW}, the simplicial group $\hat F$ is isomorphic to $F[S^1]$, the Milnor construction on the simplicial $1$-sphere which geometric realization is homotopy equivalent to $\Omega S^2$.

Let $\la\la S\ra\ra$ be the normal closure of a subset $S$ in $\hat F_{n+1}\cong F_{n}$. By the definition of the faces in formula~(\ref{equation6.6}),
\begin{equation}\label{equation6.7}
\Ker(d_i)=\la\la x_{i+1}\ra\ra
\end{equation}
in $\hat F_{n+1}\cong \pi_1(S^2\smallsetminus Q_{n+1})$ for $0\leq i\leq n$. In particular,
\begin{equation}\label{equation6.8}
\mathcal{Z}_n\hat{\mathcal{F}}=\bigcap_{i=0}^n\Ker(d_i\colon \hat F_{n+1}\to \hat F_n)=\bigcap_{i=1}^{n+1}\la\la x_i\ra\ra.
\end{equation}
By combining~\cite[Theorem 1.1]{Li-Wu} and~\cite[Theorem 4.12]{Wu1} together,
\begin{equation}\label{equation6.9}
\mathcal{B}_n\hat{\mathcal{F}}=[[\Ker(d_0),\Ker(d_1)],\ldots,\Ker(d_n)]_S.
\end{equation}

Now we start to construct an element in $\bigcap_{i=1}^{n+1}\la\la x_i\ra\ra$ that generates the quotient group
$$
\calA_{n+1}=\left.\bigcap_{i=1}^{n+1}\la\la x_i\ra\ra\right/[[\la\la x_1\ra\ra,\la\la x_2\ra\ra],\ldots,\la\la x_{n+1}\ra\ra]_S\cong \pi_n(\hat{\mathcal{F}})\cong\pi_n(\Omega S^2)
$$
for $n\leq 4$. Note that $(\hat{\mathcal{F}})_0=\hat F_1=\{1\}$ and $(\hat{\mathcal{F}})_1=\hat F_2\cong F_1=\Z$. Thus $\pi_1(\hat{\mathcal{F}})=\Z$ is generated by $x_1$. By applying Lemma~\ref{lemma6.1} together with the degeneracies in formula~(\ref{equation6.6}), we have the following:

\begin{enumerate}
\item In $\hat F_3$ with $n+1=3$, the group $\calA_3\cong \pi_2(\Omega S^2)\cong \pi_3(S^2)\cong \Z$ is generated by
$$
[s_0x_1,s_1x_1]=[x_1x_2,x_1].
$$
\item In $\hat F_4$ with $n+1=4$, the group  $\calA_4\cong \pi_3(\Omega S^2)\cong \pi_4(S^2)\cong \Z/2$ is generated by
$$
[s_0[x_1x_2,x_1],s_1[x_1x_2,x_1]]=[[x_1x_2x_3,x_2],[x_1x_2x_3,x_1]].
$$
\item In $\hat F_5$ with $n+1=5$, the group  $\calA_5\cong \pi_4(\Omega S^2)\cong \pi_5(S^2)\cong \Z/2$ is generated by
$$
\begin{array}{rl}
&[s_0[[x_1x_2x_3,x_2],[x_1x_2x_3,x_1]], s_1[[x_1x_2x_3,x_2],[x_1x_2x_3,x_1]]]\\
=&[[[x_1x_2x_3x_4,x_3],[x_1x_2x_3x_4,x_2],[[x_1x_2x_3x_4,x_2x_3],[x_1x_2x_3x_4,x_1]]].\\
\end{array}
$$
\end{enumerate}

Let $\alpha_i\in G(L_n)$ be the $i\,$th meridian such that $$p_*(\alpha_i)=x_i$$ for the epimorphism $p_*$ given in diagram~(\ref{equation5.8}). Note that $A_i=\la\la\alpha_i\ra\ra$. By equation~(\ref{equation5.9}),
$$
p_*|_{A_i}\colon A_i=\la\la \alpha_i\ra\ra\longrightarrow \la\la x_i\ra\ra
$$
is an isomorphism. The above computations give the following.
\begin{prop}
The group $\calA(L_4,L_4)=\Z/2$ is generated by
$$
l=[[\alpha_1\alpha_2\alpha_3,\alpha_2],[\alpha_1\alpha_2\alpha_3,\alpha_1]]
$$
with the $5$-link $L_4\cup l$ labeled by $\eta$ and the group $\calA(L_5,L_5)=\Z/2$ is generated by
$$
l'=[[[\alpha_1\alpha_2\alpha_3\alpha_4,\alpha_3],[\alpha_1\alpha_2\alpha_3\alpha_4,\alpha_2], [[\alpha_1\alpha_2\alpha_3\alpha_4,\alpha_2\alpha_3],[\alpha_1\alpha_2\alpha_3\alpha_4,\alpha_1]]]
$$
with the $6$-link $L_5\cup l'$ labeled by $\eta^2$.\hfill $\Box$
\end{prop}

\subsection{Some Remarks}
Let $(L,L_0)$ be any strongly nonsplittable in $S^3$ with $L\smallsetminus L_0=\{l_1,\ldots,l_n\}$. In this subsection, we provide a method how to construct the elements in $\pi_n(S^3\smallsetminus L_0)$ from $\mathcal{A}(L,L\smallsetminus L_0)$.

Let $\alpha\in \calA(L,L\smallsetminus L_0)$. Choose any element $[f]\in A(L,L\smallsetminus L_0)$ that maps to $\alpha$ in $\mathcal{A}(L,L\smallsetminus L_0)$, where
$$
f\colon S^1\longrightarrow S^3\smallsetminus L
$$
is a loop. We are going to construct certain map $S^n\to S^3\smallsetminus L_0$ from the map $f$.

\begin{lem}\label{lemma6.3}
Let $\mathbf{S}(n)=\{S(n)_{\epsilon}\}$ be a cofibrant $n$-corner with
$$
S(n)_{(0,0,\ldots,0)}\simeq S^1
$$
and $S(n)_{(\epsilon_1,\ldots,\epsilon_n)}$ is contractible for $(\epsilon_1,\ldots,\epsilon_n)\not=(0,\ldots,0)$. Then
$$
\hocolim \mathbf{S}(n)\simeq S^n
$$
for $n\geq 1$.
\end{lem}
\begin{proof}
The assertion is obvious for $n=1$. For $n=2$, the assertion follows from the homotopy pushout diagram
\begin{diagram}
S^1&\rTo&\ast\\
\dTo&&\dTo\\
\ast&\rTo&S^2.\\
\end{diagram}
For $n>2$, the assertion follows inductively by Lemma~\ref{lemma4.1}.
\end{proof}

Let $M=S^3\smallsetminus L_0$, $M_0=M_{(0,\ldots,0)}=S^3\smallsetminus L$, $M_i=M_{(0,\ldots,0,\stackrel{i}{1},0,\ldots,0)}=S^3\smallsetminus d_iL$ and
$$
M_{(\epsilon_1,\ldots,\epsilon_n)}=\bigcup_{\epsilon_i=1}M_i.
$$
Then we have the $n$-corner $\mathbf{M}(L)$ induced by the partition $(M;M_1,\ldots,M_n;M_0)$. By Theorem~\ref{theorem3.1}(ii), there is natural isomorphism
\begin{equation}\label{equation6.10}
\rho_{\mathbf{M}(L)}\colon \calA(L,L\smallsetminus L_0)\longrightarrow \pi_n(M).
\end{equation}

\begin{lem}\label{lemma6.4}
Let $f\colon S^1\to S^3\smallsetminus L$ be a loop such that $[f]\in A(L,L\smallsetminus L_0)$ and let $\mathbf{S}(n)=\{S(n)_{\epsilon}\}$ be a cofibrant $n$-corner with
$$
S(n)_{(0,0,\ldots,0)}= S^1
$$
and $S(n)_{(\epsilon_1,\ldots,\epsilon_n)}$ is contractible for $(\epsilon_1,\ldots,\epsilon_n)\not=(0,\ldots,0)$. Suppose that $(L,L_0)$ is strongly nonsplittable. Then there exists a morphism of $n$-corners
$$
\mathbf{f}\colon \mathbf{S}(n)\longrightarrow \mathbf{M}(L)
$$
such that
$$
f_{(0,\ldots,0)}=f\colon S(n)_{(0,0,\ldots,0)}= S^1\longrightarrow M_{(0,\ldots,0)}=S^3\smallsetminus L.
$$
\end{lem}
\begin{proof}
The proof is given inductively by constructing certain extension maps.
Since $[f]\in A(L,L\smallsetminus L_0)$, the composite
$$
S(n)_{(0,0,\ldots,0)}=S^1\rTo^{f=f_{(0,\ldots,0)}} M_{(0,\ldots,0)}=S^3\smallsetminus L\rInto S^3\smallsetminus d_iL
$$
is null homotopic for each $1\leq i\leq n$. Thus there exists a map $f_{(0,\ldots,0,\stackrel{i}{1},0,\ldots,0)}$ such that the diagram
\begin{diagram}
S(n)_{(0,0,\ldots,0)}&\rTo^{f_{(0,\ldots,0)}}& M_{(0,\ldots,0)}\\
\dInto&&\dInto\\
S(n)_{(0,\ldots,0,\stackrel{i}{1},0,\ldots,0)} &\rDashto^{f_{(0,\ldots,0,\stackrel{i}{1},0,\ldots,0)}}& M_{(0,\ldots,0,\stackrel{i}{1},0,\ldots,0)}\\
\end{diagram}
commutes (strictly) for each $1\leq i\leq n$. Now suppose that
$$
f_{\alpha}\colon S(n)_{\alpha}\longrightarrow M_{\alpha}
$$
has been constructed for $\alpha<\epsilon\in \{0,1\}^n\smallsetminus\{(1,\ldots,1)\}$. Let
$$\calC_{\epsilon}=\{\alpha \in \{0,1\}^n\ |\ \alpha<\epsilon\}$$ be the subcategory of $\{0,1\}^n$ consisting of those object less than $\epsilon$.
If $\calC$ has no non-identity morphism, then $\epsilon=(0,\ldots,0,\stackrel{i}{1},0,\ldots,0)$ for some $1\leq i\leq n$. In this case, we have constructed the maps $f_{(0,\ldots,0,\stackrel{i}{1},0,\ldots,0)}$. Thus we may assume that $\calC_{\epsilon}$ has at least one non-identity morphism. Let $\epsilon=(\epsilon_1,\ldots,\epsilon_n)$. Observe that
$$
\calC_{\epsilon}=\{\alpha=(\alpha_1,\ldots,\alpha_n) \in \{0,1\}^n \ |\ \alpha_i\leq \epsilon_i \textrm{ for } 1\leq i\leq n \textrm{ with } \alpha\not=\epsilon\}.
$$
The diagrams
$
\mathcal{S}(n)|_{\calC_{\epsilon}}\textrm{ and } \mathcal{M}(L)|_{\calC_{\epsilon}}
$
are $t$-corners with $t\geq 2$. From the induction hypothesis, $f_{\alpha}\colon S(n)_{\alpha}\longrightarrow M_{\alpha}$ has been constructed and so it induces a map
$$
\colim f_{\alpha}\colon \colim_{\calC_{\epsilon}}\mathcal{S}(n)|_{\calC_{\epsilon}}\longrightarrow \colim_{\calC_{\epsilon}} \mathcal{M}(L)|_{\calC_{\epsilon}}=M_{\epsilon},
$$
where $\colim_{\calC_{\epsilon}} \mathcal{M}(L)|_{\calC_{\epsilon}}=M_{\epsilon}$ because $\mathcal{M}(L)$ is induced by the partition. By Lemma~\ref{lemma6.3},
$$
\colim_{\calC_{\epsilon}}\mathcal{S}(n)|_{\calC_{\epsilon}}\simeq S^t.
$$
Since $t\geq 2$, we have $\pi_t(M_{\epsilon})=0$ by the assumption of strongly nonsplittability. It follows that there exists an extension
\begin{diagram}
\colim_{\calC_{\epsilon}}\mathcal{S}(n)|_{\calC_{\epsilon}}&\rTo^{\colim f_{\alpha}}& M_{\epsilon}\\
\dInto&\ruDashto>{f_{\epsilon}}&\\
S(n)_{\epsilon}.&&\\
\end{diagram}
The induction is finished and hence the result.
\end{proof}
Now we give the following construction:
\begin{enumerate}
\item[] \textit{Let $f\colon S^1\to S^3\smallsetminus L$ be a loop such that $[f]\in A(L,L\smallsetminus L_0)$ and let
$$
\mathbf{f}\colon \mathbf{S}(n)\longrightarrow \mathbf{M}(L)
$$
be any morphism of $n$-corners such that
$$
f_{(0,\ldots,0)}=f\colon S(n)_{(0,0,\ldots,0)}= S^1\longrightarrow M_{(0,\ldots,0)}=S^3\smallsetminus L.
$$
Define
$$
\theta(f)=\colim\mathbf{f}\colon \colim\mathbf{S}(n)\simeq\hocolim \mathbf{S}(n)\simeq S^n\longrightarrow \colim \mathbf{M}(L)=M=S^3\smallsetminus L_0.
$$}
\end{enumerate}
The map $\theta(f)$ is of course dependent on the choice of the morphism $\mathbf{f}$. But its homotopy class is independent on the choice of $\mathbf{f}$.

\begin{thm}\label{theorem6.5}
Let $(L,L_0)$ be a strongly nonsplittable pair of links in $S^3$ with $L\smallsetminus L_0$ an $n$-link. Let $f\colon S^1\to S^3\smallsetminus L$ be a loop such that $[f]\in A(L,L\smallsetminus L_0)$ that represents $\alpha$ and let $\theta(f)$ be defined as above. Then
$$
[\theta(f)]=\rho_{\mathbb{M}}(\alpha).
$$
\end{thm}
\begin{proof}
Observe that the $n$-corner $\mathcal{S}(n)$ satisfies the connectivity hypothesis in~\cite[Theorem 1]{EM} with
$$
\Ker(\pi_1(S(n)_{(0,\ldots,0)})\to \pi_1(S(n)_{(0,\ldots,0,\stackrel{i}{1},0,\ldots,0)}))=\pi_1(S(n)_{(0,\ldots,0)})=\pi_1(S^1)=\Z
$$
because $S(n)_{(\epsilon_1,\ldots,\epsilon_n)})$ is contractible for $(\epsilon_1,\ldots,\epsilon_n)\not=(0,\ldots,0), (1,\ldots,1)$. By~\cite[Theorem 1]{EM}, there is an isomorphism
$$
\rho \colon \pi_1(S^1)\longrightarrow \pi_n(\hocolim \mathcal{S}(n)).
$$
From the naturality of $\rho$, there is a commutative diagram
\begin{diagram}
\pi_1(S^1)=\Z&\rTo^{\rho_{\mathcal{S}(n)}}& \pi_n(\hocolim \mathcal{S}(n))=\pi_n(S^2)\\
\dTo>{f_{*}=f_{(0,\ldots,0)_*}}&&\dTo>{\theta(f)_*}\\
\mathcal{A}(L,L_0)&\rTo^{\rho_{\mathbf{M}(L)}}&\pi_n(S^3\smallsetminus L_0).\\
\end{diagram}
It follows that
$$
\begin{array}{rcl}
[\theta(f)]&=&\theta(f)_*(\iota_n)\\
&=&\theta(f)_*(\rho_{\mathcal{S}(n)}(\iota_1))\\
&=&\rho_{\mathbf{M}(L)}(f_*(\iota_1)) \\
&=&\rho_{\mathbf{M}(L)}(\alpha)\\
\end{array}
$$
and hence the result.
\end{proof}

\begin{rem}{\rm
The construction $\theta(f)$ is an analogue of the Massey product. By Theorems~\ref{theorem1.1} and ~\ref{theorem6.5}, all elements in $\pi_*(S^3)$ can be obtained by this kind of operations. Moreover the mapping cone of the morphism
$$
\mathbf{f}\colon \mathbf{S}(n)\longrightarrow \mathbf{M}(L)
$$
gives a cubical resolution for the $2$-cell complex given by the homotopy cofibre of $\theta(f)\colon S^n\to S^3\smallsetminus L_0$, which might be useful for studying the $2$-complexes given in the form of $(S^3\smallsetminus L_0)\cup e^{n+1}$.
}\end{rem}

\noindent\small{\textbf{Acknowledgements.}} The author would like to thank Joan Birman and Haynes Miller for their encouragements and help suggestions on this project.

\end{document}